%% file: root.tex
\documentclass[article]{IEEEtran}
\usepackage{cite}
\usepackage{color}
\usepackage[pdftex]{graphicx}
\graphicspath{{fig/}{jpeg/}}
\usepackage[cmex10]{amsmath}
\usepackage{amssymb}
\usepackage{algorithm}
\usepackage{algorithmic}
\usepackage{amsmath,amssymb,lipsum}
\usepackage{hyperref}
\usepackage{comment}
\usepackage{enumerate}
\usepackage{graphicx}
\usepackage{subcaption}
\usepackage{setspace}
\usepackage{caption}
\input{mysymbol.sty}
\input{my_sections.tex}

\renewcommand{\QED}{\hfill\ensuremath{\blacksquare}}
\usepackage{theorem}
\newtheorem{thm}{Theorem}
\newtheorem{lemma}{Lemma}
\newtheorem{proposition}{Proposition}
\theoremstyle{definition}
\newtheorem{assumption}{Assumption}
\newtheorem{corollary}{Corollary}

\title{DQM: Decentralized Quadratically Approximated Alternating Direction Method of Multipliers}

\author{Aryan Mokhtari, Wei Shi, Qing Ling, and Alejandro Ribeiro
\thanks{{Work supported by NSF CAREER CCF-0952867, ONR N00014-12-1-0997, and NSFC 61004137. A. Mokhtari and A. Ribeiro are with the Dept. of Electrical and Systems Engineering, University of Pennsylvania, 200 S 33rd St., Philadelphia, PA 19104. Email: \{aryanm, aribeiro\}@seas.upenn.edu. W. Shi and Q. Ling are with the Dept. of Automation, University of Science and Technology of China, 96 Jinzhao Rd., Hefei, Anhui, 230026, China. Email: \{shiwei00, qingling\}@mail.ustc.edu.cn. This paper expands results and presents proofs that were preliminarily reported in \cite{DQMglobalSIP}.
}}
}

\begin{document}

\maketitle
\thispagestyle{empty}

\begin{abstract}
This paper considers decentralized consensus optimization problems where nodes of a network have access to different summands of a global objective function. Nodes cooperate to minimize the global objective by exchanging information with neighbors only. A decentralized version of the alternating directions method of multipliers (DADMM) is a common method for solving this category of problems. DADMM exhibits linear convergence rate to the optimal objective but its implementation requires solving a convex optimization problem at each iteration. This can be computationally costly and may result in large overall convergence times. The decentralized quadratically approximated ADMM algorithm (DQM), which minimizes a quadratic approximation of the objective function that DADMM minimizes at each iteration, is proposed here. The consequent reduction in computational time is shown to have minimal effect on convergence properties. Convergence still proceeds at a linear rate with a guaranteed constant that is asymptotically equivalent to the DADMM linear convergence rate constant. Numerical results demonstrate advantages of DQM relative to DADMM and other alternatives in a logistic regression problem.
\end{abstract}

\begin{keywords}
Multi-agent network, decentralized optimization, Alternating Direction Method of Multipliers.
\end{keywords}

\input{Introduction.tex}

\input{Problem.tex}
\input{Con_Anlysis.tex}
\input{Simulations.tex}
\input{Conclusions.tex}
\input{Appendix.tex}
\bibliographystyle{IEEEtran}
  \bibliography{bmc_article}
   \end{document}

%% file: my_sections.tex
\usepackage{needspace}

\newcommand{\myparagraph}[1]{\needspace{1\baselineskip}\medskip\noindent {\it #1.}}



%% file: Introduction.tex

\section{Introduction}\label{sec_Introduction}

Decentralized algorithms are used to solve optimization problems
where components of the objective are available at different nodes
of a network. Nodes access their local cost functions only but
try to minimize the aggregate cost by exchanging information with
their neighbors. Specifically, consider a variable
$\tbx\in\reals^p$ and a connected network containing $n$ nodes
each of which has access to a local cost function
$f_i:\reals^p\to\reals$. The nodes' goal is to find the optimal
argument of the global cost function $ \sum_{i=1}^{n}f_i(\tbx)$,
\begin{equation}\label{original_optimization_problem1}
  \tbx^*
      \  =\ \argmin_{\tbx}\sum_{i=1}^{n} f_i(\tbx).
\end{equation}
Problems of this form arise in, e.g., decentralized control
\cite{Bullo2009,Cao2013-TII,LopesEtal8}, wireless communication
\cite{Ribeiro10,Ribeiro12}, sensor networks
\cite{Schizas2008-1,KhanEtal10,cRabbatNowak04}, and large scale
machine learning
\cite{bekkerman2011scaling,Tsianos2012-allerton-consensus,Cevher2014}.
In this paper we assume that the local costs $f_i$ are twice
differentiable and strongly convex.

There are different algorithms to solve
\eqref{original_optimization_problem1} in a decentralized manner
which can be divided into two major categories. The ones that
operate in the primal domain and the ones that operate in the dual
domain. Among primal domain algorithms, decentralized
(sub)gradient descent (DGD) methods are well studied
\cite{Nedic2009,Jakovetic2014-1,YuanQing}. They can be interpreted
as either a mix of local gradient descent steps with successive
averaging or as a penalized version of
\eqref{original_optimization_problem1} with a penalty term that
encourages agreement between adjacent nodes. This latter
interpretation has been exploited to develop the network Newton
(NN) methods that attempt to approximate the Newton step of this
penalized objective in a distributed manner
\cite{NN-part1,NN-part2}. The methods that operate in the dual
domain consider a constraint that enforces equality between nodes'
variables. They then ascend on the dual function to find optimal
Lagrange multipliers with the solution of
\eqref{original_optimization_problem1} obtained as a byproduct
\cite{Schizas2008-1,BoydEtalADMM11,Shi2014-ADMM,rabbat2005generalized}.
Among dual descent methods, decentralized implementation of the
alternating directions method of multipliers (ADMM), known as
DADMM, is proven to be very efficient with respect to convergence
time \cite{Schizas2008-1,BoydEtalADMM11,Shi2014-ADMM}.

A fundamental distinction between primal methods such as DGD and
NN and dual domain methods such as DADMM is that the former
compute local gradients and Hessians at each iteration while the
latter minimize local pieces of the Lagrangian at each step --
this is necessary since the gradient of the dual function is
determined by Lagrangian minimizers. Thus, iterations in dual
domain methods are, in general, more costly because they require
solution of a convex optimization problem. However, dual methods
also converge in a smaller number of iterations because they
compute approximations to $\tbx^*$ instead of descending towards
$\tbx^*$. Having complementary advantages, the choice between
primal and dual methods depends on the relative cost of
computation and communication for specific problems and platforms.
Alternatively, one can think of developing methods that combine
the advantages of ascending in the dual domain without requiring
solution of an optimization problem at each iteration. This can be
accomplished by the decentralized linearized ADMM (DLM) algorithm
\cite{cQingRibeiroADMM14,ling2014dlm}, which replaces the
minimization of a convex objective required by ADMM with the
minimization of a first order linear approximation of the
objective function. This yields per-iteration problems that can be
solved with a computational cost akin to the computation of a
gradient and a method with convergence properties closer to DADMM
than DGD.

If a first order approximation of the objective is useful, a
second order approximation should decrease convergence times
further. The decentralized quadratically approximated ADMM (DQM)
algorithm that we propose here minimizes a quadratic approximation
of the Lagrangian minimization of each ADMM step. This quadratic
approximation requires computation of local Hessians but results
in an algorithm with convergence properties that are: (i) better than the convergence properties of DLM; (ii) asymptotically
identical to the convergence behavior of DADMM. The technical
contribution of this paper is to prove that (i) and (ii) are true
from both analytical and practical perspectives.

We begin the paper by discussing solution of
\eqref{original_optimization_problem1} with DADMM and its
linearized version DLM (Section \ref{sec_preliminaries}). Both of
these algorithms perform updates on dual and primal auxiliary
variables that are identical and computationally simple. They
differ in the manner in which principal primary variables are
updated. DADMM solves a convex optimization problem and DLM solves
a regularized linear approximation. We follow with an explanation
of DQM that differs from DADMM and DLM in that it minimizes a
quadratic approximation of the convex problem that DADMM solves
exactly and DLM approximates linearly (Section \ref{sec:DQM}). We
also explain how DQM can be implemented in a distributed manner
(Proposition \ref{update_system_prop} and Algorithm
\ref{algo_DQM}). Convergence properties of DQM are then analyzed
(Section \ref{sec:Analysis}) where linear convergence is
established (Theorem \ref{DQM_convergence} and Corollary
\ref{convergence_of_primal}). Key in the analysis is the error
incurred when approximating the exact minimization of DADMM with
the quadratic approximation of DQM. This error is shown to
decrease as iterations progress (Proposition
\ref{error_vector_proposition}) faster than the rate that the error of DLM approaches zero (Proposition
\ref{error_vector_proposition_2}). This results in DQM having a
guaranteed convergence constant strictly smaller than the DLM
constant that approaches the guaranteed constant of DADMM for
large iteration index (Section \ref{sec:rate_comparison}). We
corroborate analytical results with numerical evaluations in a
logistic regression problem (Section \ref{sec:simulations}). We
show that DQM does outperform DLM and show that convergence paths
of DQM and DADMM are almost identical (Section \ref{comparison}).
Overall computational cost of DQM is shown to be smaller, as
expected.

\myparagraph{\bf Notation} Vectors are written as
$\bbx\in\reals^n$ and matrices as $\bbA\in\reals^{n\times n}$.
Given $n$ vectors $\bbx_i$, the vector
$\bbx=[\bbx_1;\ldots;\bbx_n]$ represents a stacking of the
elements of each individual $\bbx_i$. We use $\|\bbx\|$ to denote
the Euclidean norm of vector $\bbx$ and $\|\bbA\|$ to denote the
Euclidean norm of matrix $\bbA$. The gradient of a function $f$ at
point $\bbx$ is denoted as $\nabla f(\bbx)$ and the Hessian is
denoted as $\nabla^2 f(\bbx)$. We use $\sigma(\bbB)$ to denote the singular values of matrix $\bbB$ and $\lambda(\bbA)$ to denote the eigenvalues of matrix $\bbA$.

%% file: Problem.tex

%
\section{Distributed Alternating Directions Method of Multipliers} \label{sec_preliminaries}

Consider a connected network with $n$ nodes and $m$ edges where the set of nodes is $\mathcal{V}=\{1, \dots, n\}$ and the set of ordered edges $\mathcal{E}$ contains pairs $(i,j)$ indicating that $i$ can communicate to $j$. We restrict attention to symmetric networks in which $(i,j)\in\mathcal{E}$ if and only if $(j,i)\in\mathcal{E}$ and define node $i$'s neighborhood as the set $\mathcal{N}_i=\{j\mid (i,j)\in \mathcal{E}\}$. In problem \eqref{original_optimization_problem1} agent $i$ has access to the local objective function $f_{i}(\tbx)$ and agents cooperate to minimize the global cost $\sum_{i=1}^nf_{i}(\tbx)$. This specification is more naturally formulated by defining variables $\bbx_i$ representing the local copies of the variable $\tbx$. We also define the auxiliary variables $\bbz_{ij}$ associated with edge $(i,j)\in \mathcal{E}$ and rewrite \eqref{original_optimization_problem1} as
\begin{alignat}{2}\label{original_optimization_problem2}
   \{\bbx_i^*\}_{i=1}^n :=\
       &\argmin_{\bbx} &&\ \sum_{i=1}^{n}\ f_{i}(\bbx_{i}),         \\ \nonumber
       &\st            &&\ \bbx_{i}=\bbz_{ij},\ \bbx_{j}=\bbz_{ij}, \
                           \text{for all\ } (i, j)\in\ccalE .
\end{alignat}
The constraints $\bbx_{i}=\bbz_{ij}$ and $ \bbx_{j}=\bbz_{ij}$ enforce that the variable $\bbx_i$ of each node $i$ is equal to the variables $\bbx_j$ of its neighbors $j\in \ccalN_i$. This condition in association with network connectivity implies that a set of variables $\{\bbx_1,\dots,\bbx_n\}$ is feasible for problem \eqref{original_optimization_problem2} if and only if all the variables $\bbx_i$ are equal to each other, i.e., if $\bbx_1=\dots=\bbx_n$. Therefore, problems \eqref{original_optimization_problem1} and \eqref{original_optimization_problem2} are equivalent in the sense that for all $i$ and $j$ the optimal arguments of \eqref{original_optimization_problem2} satisfy $\bbx_i^*=\tbx^*$ and $\bbz_{ij}=\tbx^*$, where $\tbx^*$ is the optimal argument of \eqref{original_optimization_problem1}.

To write problem \eqref{original_optimization_problem2} in a
matrix form, define $\bbA_s\in \reals^{mp\times np} $ as the block
source matrix which contains $m\times n$ square blocks
$(\bbA_s)_{e,i}\in \reals^{p\times p}$. The block $(\bbA_s)_{e,i}$
is not identically null if and only if the edge $e$ corresponds to
$e=(i,j)\in \ccalE$ in which case $(\bbA_s)_{e,i}=\bbI_{p}$.
Likewise, the block destination matrix $\bbA_d\in \reals^{mp\times
np} $ contains  $m\times n$ square blocks $(\bbA_d)_{e,i}\in
\reals^{p\times p}$. The square block $(\bbA_d)_{e,i}=\bbI_p$ when
$e$ corresponds to $e=(j,i)\in \ccalE$ and is null otherwise.
Further define $\bbx:=[\bbx_1;\dots;\bbx_n]\in \reals^{np}$ as a
vector concatenating all local variables $\bbx_i$, the vector
$\bbz:=[\bbz_1;\dots;\bbz_m]\in \reals^{mp}$ concatenating all
auxiliary variables $\bbz_e=\bbz_{ij}$, and the  aggregate
function $f:\reals^{np}\to\reals$ as $f(\bbx):=
\sum_{i=1}^{n}f_i(\bbx_i)$. We can then rewrite
\eqref{original_optimization_problem2} as 
\begin{equation}\label{original_optimization_problem3}
   \bbx^* := \argmin_{\bbx}f(\bbx), \quad
             \st\ \bbA_s\bbx-\bbz=\bb0, \
                  \bbA_d\bbx-\bbz=\bb0.
\end{equation}
Define now the matrix $\bbA=[\bbA_s;\bbA_d]\in \reals^{2mp\times np}$ which stacks the source and destination matrices, and the matrix $\bbB=[-\bbI_{mp};-\bbI_{mp}]\in \reals^{2mp\times mp}$ which stacks two negative identity matrices of size $mp$ to rewrite \eqref{original_optimization_problem3} as
\begin{align}\label{original_optimization_problem4}
   \bbx^* := \argmin_{\bbx}f(\bbx), \quad
             \st\ \bbA\bbx+\bbB\bbz=\bb0.
\end{align}
DADMM is the application of ADMM to solve \eqref{original_optimization_problem4}. To develop this algorithm introduce Lagrange multipliers $\bbalpha_e=\bbalpha_{ij}$ and $\bbbeta_e=\bbbeta_{ij}$ associated with the constraints $\bbx_{i}=\bbz_{ij}$ and $\bbx_{j}=\bbz_{ij}$ in \eqref{original_optimization_problem2}, respectively. Define $\bbalpha:=[\bbalpha_1;\dots;\bbalpha_m]$ as the concatenation of the multipliers $\bbalpha_e$ which yields the multiplier of the constraint $\bbA_s\bbx-\bbz=\bb0$ in \eqref{original_optimization_problem3}. Likewise, the corresponding Lagrange multiplier of the constraint $\bbA_d\bbx-\bbz=\bb0$ in \eqref{original_optimization_problem3} can be obtained by stacking the multipliers $\bbbeta_e$ to define $\bbbeta:=[\bbbeta_1;\dots;\bbbeta_m]$. Grouping $\bbalpha$ and $\bbbeta$ into $\bblambda:=[\bbalpha;\bbbeta]\in \reals^{2mp}$ leads to the Lagrange multiplier $\bblambda$ associated with the constraint $\bbA\bbx+\bbB\bbz=\bb0$ in \eqref{original_optimization_problem4}. Using these definitions and introducing a positive constant $c>0$ we write the augmented Lagrangian of \eqref{original_optimization_problem4} as
\begin{equation}\label{lagrangian}
\ccalL(\bbx,\bbz,\bblambda) := f(\bbx)+\bblambda^{T}\left(\bbA\bbx+\bbB\bbz\right)+ \frac{c}{2}\left\|\bbA\bbx+\bbB\bbz\right\|^2.
\end{equation}
The idea of ADMM is to minimize the Lagrangian $\ccalL(\bbx,\bbz,\bblambda)$ with respect to $\bbx$, follow by minimizing the updated Lagrangian with respect to $\bbz$, and finish each iteration with an update of the multiplier $\bblambda$ using dual ascent. To be more precise, consider the time index $k \in \naturals$ and define $\bbx_k$, $\bbz_k$, and $\bblambda_k$ as the iterates at step $k$. At this step, the augmented Lagrangian is minimized with respect to $\bbx$ to obtain the iterate
\begin{equation}\label{ADMM_x_update}
\bbx_{k+1}= \argmin_{\bbx}  f(\bbx)+\bblambda_k^{T}\left(\bbA\bbx+\bbB\bbz_k\right)+ \frac{c}{2}\left\|\bbA\bbx+\bbB\bbz_k\right\|^2.
\end{equation}
Then, the augmented Lagrangian is minimized with respect to the auxiliary variable $\bbz$ using the updated variable $\bbx_{k+1}$ to obtain
\begin{align}\label{ADMM_z_update}
\bbz_{k+1}= \argmin_{\bbz}  &\ f(\bbx_{k+1}) \\ \nonumber
                &+\bblambda_k^{T}\left(\bbA\bbx_{k+1}+\bbB\bbz\right)
                + \frac{c}{2}\left\|\bbA\bbx_{k+1}+\bbB\bbz\right\|^2 .
\end{align}
After updating the variables $\bbx$ and $\bbz$, the Lagrange multiplier $\bblambda_k$ is updated through the dual ascent iteration
\begin{equation}\label{ADMM_lambda_update}
\bblambda_{k+1}=\bblambda_{k}+c\left(\bbA\bbx_{k+1}+\bbB\bbz_{k+1}\right).
\end{equation}
The DADMM algorithm is obtained by observing that the structure of the matrices $\bbA$ and $\bbB$ is such that \eqref{ADMM_x_update}-\eqref{ADMM_lambda_update} can be implemented in a distributed manner \cite{Schizas2008-1,BoydEtalADMM11,Shi2014-ADMM}.

The updates for the auxiliary variable $\bbz$ and the Lagrange multiplier $\bblambda$ are not costly in terms of computation time. However, updating the primal variable $\bbx$ can be expensive as it entails the solution of an optimization problem [cf. \eqref{ADMM_x_update}]. The DLM algorithm avoids this cost with an inexact update of the primal variable iterate $\bbx_{k+1}$. This inexact update relies on approximating the aggregate function value $f(\bbx_{k+1})$ in \eqref{ADMM_x_update} through a regularized linearization of the aggregate function $f$ in a neighborhood of the current variable $\bbx_{k}$. This regularized approximation takes the form $f(\bbx)\approx f(\bbx_k)+\nabla f(\bbx_k)^T(\bbx-\bbx_k)+(\rho/2)\|\bbx-\bbx_k\|^2$ for a given positive constant $\rho>0$. Consequently, the update formula for the primal variable $\bbx$ in DLM replaces the DADMM exact minimization in \eqref{ADMM_x_update} by the minimization of the quadratic form
\begin{align}\label{DLM_x_update}
\bbx_{k+1}= \argmin_{\bbx}\ & f(\bbx_k)+\nabla f(\bbx_k)^T(\bbx-\bbx_k)+\frac{\rho}{2}\|\bbx-\bbx_k\|^2
\nonumber \\
& \quad
+\bblambda_k^{T}\left(\bbA\bbx+\bbB\bbz_k\right)+ \frac{c}{2}\left\|\bbA\bbx+\bbB\bbz_k\right\|^2.
\end{align}
The first order optimality condition for \eqref{DLM_x_update} implies that the updated variable $\bbx_{k+1}$ satisfies
\begin{equation}\label{DLM_optimality_cond}
\nabla f(\bbx_{k})+\rho (\bbx_{k+1}-\bbx_{k})+\bbA^T\bblambda_{k}+c \bbA^{T} \left( \bbA\bbx_{k+1}+\bbB \bbz_{k}  \right)=\bb0.
\end{equation}
According to \eqref{DLM_optimality_cond}, the updated variable $\bbx_{k+1}$ can be computed by inverting the positive definite matrix $\rho\bbI+c\bbA^T\bbA$. This update can also be implemented in a distributed manner.

The sequence of variables $\bbx_k$ generated by DLM converges linearly to the optimal argument $\bbx^*$ \cite{cQingRibeiroADMM14}. Although this is the same rate of DADMM, linear convergence constant of DLM is smaller than the one for DADMM (see Section \ref{sec:rate_comparison}), and can be much smaller depending on the condition number of the local functions $f_i$ (see Section \ref{comparison}). To close the gap between these constants we can use a second order approximation of \eqref{ADMM_x_update}. This is the idea of DQM that we introduce in the following section.

\section{DQM: Decentralized Quadratically Approximated ADMM}\label{sec:DQM}

DQM uses a local quadratic approximation of the primal function $f(\bbx)$ around the current iterate $\bbx_k$. If we let $\bbH_{k}:=\nabla^2 f(\bbx_{k})$ denote the primal function Hessian evaluated at $\bbx_k$ the quadratic approximation of $f$ at $\bbx_k$ is  $f(\bbx) \approx f(\bbx_k)+\nabla f(\bbx_k)^T(\bbx-\bbx_k)+(1/2)(\bbx-\bbx_k)^T\bbH_k(\bbx-\bbx_k)$. Using this approximation in \eqref{ADMM_x_update} yields the DQM update that we therefore define as
\begin{align}\label{DQM_x_update}
    \bbx_{k+1}:= \argmin_{\bbx} f(&\bbx_k)
                    +\nabla f(\bbx_k)^T(\bbx-\bbx_k)                  \\ \nonumber
                  & +\frac{1}{2}(\bbx-\bbx_k)^T\bbH_k(\bbx-\bbx_k)    \\ \nonumber
                  & +\bblambda_k^{T}\left(\bbA\bbx+\bbB\bbz_k\right)
                    + \frac{c}{2}\left\|\bbA\bbx+\bbB\bbz_k\right\|^2 .
\end{align}
Comparison of \eqref{DLM_x_update} and \eqref{DQM_x_update} shows
that in DLM the quadratic term $(\rho/2)\|\bbx_{k+1}-\bbx_{k}\|^2$
is added to the first-order approximation of the primal objective
function, while in DQM the second order approximation of the
primal objective function is used to reach a more accurate
approximation for $f(\bbx)$. Since \eqref{DQM_x_update} is a
quadratic program, the first order optimality condition yields a
system of linear equations that can be solved to find
$\bbx_{k+1}$,
\begin{align}\label{DQM_update_1}
    \nabla f(\bbx_{k})+ \bbH_{k}(\bbx_{k+1}\!-\!\bbx_{k})
   +\bbA^T\!\bblambda_{k}+c \bbA^{T}\! \left( \bbA\bbx_{k+1}
   +\bbB \bbz_{k}  \right)=\bb0.
\end{align}
This update can be solved by inverting the matrix $\bbH_k+c\bbA^T\bbA$ which is invertible if, as we are assuming, $f(\bbx)$ is strongly convex.

The DADMM updates in \eqref{ADMM_z_update} and \eqref{ADMM_lambda_update} are used verbatim in DQM, which is therefore defined by recursive application of \eqref{DQM_update_1}, \eqref{ADMM_z_update}, and \eqref{ADMM_lambda_update}. It is customary to consider the first order optimality conditions of \eqref{ADMM_z_update} and to reorder terms in \eqref{ADMM_lambda_update} to rewrite the respective updates as
\begin{align}\label{DQM_update_2}
   \bbB^T\bblambda_{k} + c \bbB^T\left( \bbA\bbx_{k+1}+\bbB\bbz_{k+1} \right)  &=\bbzero,
   \nonumber\\
   \bblambda_{k+1}-\bblambda_{k}-c\left( \bbA\bbx_{k+1}+\bbB\bbz_{k+1} \right) &=\bbzero.
\end{align}
DQM is then equivalently defined by recursive solution of the system of linear equations in \eqref{DQM_update_1} and \eqref{DQM_update_2}. This system, as is the case of DADMM and DLM, can be reworked into a simpler form that reduces communication cost. To derive this simpler form we assume a specific structure for the initial vectors $\bblambda_{0}=[\bbalpha_0;\bbbeta_0]$, $\bbx_0$, and $\bbz_0$ as introduced in the following assumption.

%
\begin{assumption}\label{initial_val_assum}
Define the oriented incidence matrix as $\bbE_o:=\bbA_s-\bbA_d$ and the unoriented incidence matrix as $\bbE_u:=\bbA_s+\bbA_d$. The initial Lagrange multipliers $\bbalpha_0$ and $\bbbeta_0$, and the initial variables $\bbx_0$ and $\bbz_0$ are chosen such that:
\begin{enumerate}[(a)]
\item The multipliers are opposites of each other, $\bbalpha_0=-\bbbeta_0$.
\item The initial primal variables satisfy $\bbE_u\bbx_0=2\bbz_0$.
\item The initial multiplier $\bbalpha_0$ lies in the column space of $\bbE_o$.
\end{enumerate} \end{assumption}

%
Assumption \ref{initial_val_assum} is minimally restrictive. The only non-elementary condition is (c) but that can be satisfied by $\bbalpha_0=\bbzero$. Nulling all other variables, i.e., making $\bbbeta_0=\bbzero$, $\bbx_0=\bbzero$, and $\bbz_0=\bbzero$ is a trivial choice to comply with conditions (a) and (b) as well. An important consequence of the initialization choice in \eqref{initial_val_assum} is that if the conditions in Assumption \ref{initial_val_assum} are true at time $k=0$ they stay true for all subsequent iterations $k>0$ as we state next.

%
\begin{lemma}\label{lag_var_lemma}
Consider the DQM algorithm as defined by
\eqref{DQM_update_1}-\eqref{DQM_update_2}. If Assumption
\ref{initial_val_assum} holds, then for all $k\geq 0$ the Lagrange
multipliers $\bbalpha_k$ and $\bbbeta_k$, and the variables
$\bbx_k$ and $\bbz_k$ satisfy:
\begin{enumerate}[(a)]
\item The multipliers are opposites of each other, $\bbalpha_k=-\bbbeta_k$.
\item The primal variables satisfy $\bbE_u\bbx_k=2\bbz_k$.
\item The multiplier $\bbalpha_k$ lies in the column space of $\bbE_o$.
\end{enumerate}
\end{lemma}

%
\begin{myproof} See Appendix \ref{app_lag_var_lemma}. \end{myproof}

%
The validity of (c) in Lemma \ref{lag_var_lemma} is important for the convergence analysis of Section \ref{sec:Analysis}. The validity of (a) and (b) means that maintaining multipliers $\bbalpha_k$ and $\bbbeta_k$ is redundant because they are opposites and that maintaining variables $\bbz_k$ is also redundant because they can be computed as $\bbz_k = \bbE_u\bbx_k/2$. It is then possible to replace \eqref{DQM_update_1}-\eqref{DQM_update_2} by a simpler system of linear equations as we explain in the following proposition.

%
\begin{proposition}\label{update_system_prop}
Consider the DQM algorithm as defined by \eqref{DQM_update_1}-\eqref{DQM_update_2} and define the sequence $\bbphi_{k}:=\bbE_{o}^T\bbalpha_k$. Further define the unoriented Laplacian as $\bbL_u:=(1/2)\bbE_u^T\bbE_u$, the oriented Laplacian as $\bbL_o=(1/2)\bbE_o^T\bbE_o$, and the degree matrix as $\bbD:=(\bbL_u+\bbL_o)/2$. If Assumption \ref{initial_val_assum} holds true, the DQM iterates $\bbx_k$ can be generated as
\begin{align}\label{x_update_formula}
   \bbx_{k+1}   &= (2c\bbD+\bbH_k)^{-1}
                   \left[(c\bbL_u+\bbH_k)\bbx_k
                           -\nabla f(\bbx_{k})-\bbphi_k \right], \nonumber\\
   \bbphi_{k+1} &= \bbphi_k +c \bbL_o\bbx_{k+1}.
\end{align}
\end{proposition}

%
\begin{myproof} See Appendix \ref{app_update_system}. \end{myproof}

%
Proposition \ref{update_system_prop} states that by introducing the sequence of variables $\bbphi_{k}$, the DQM primal iterates $\bbx_k$ can be computed through the recursive expressions in \eqref{x_update_formula}. These recursions are simpler than \eqref{DQM_update_1}-\eqref{DQM_update_2} because they eliminate the auxiliary variables $\bbz_k$ and reduce the dimensionality of $\bblambda_k$ -- twice the number of edges -- to that of $\bbphi_k$ -- the number of nodes. Further observe that if \eqref{x_update_formula} is used for implementation we don't have to make sure that the conditions of Assumption \ref{initial_val_assum} are satisfied. We just need to pick $\bbphi_{0}:=\bbE_{o}^T\bbalpha_0$ for some $\bbalpha_0$ in the column space of $\bbE_0$ -- which is not difficult, we can use, e.g., $\bbphi_0=\bbzero$. The role of Assumption \ref{initial_val_assum} is to state conditions for which the expressions in \eqref{DQM_update_1}-\eqref{DQM_update_2} are an equivalent representation of \eqref{x_update_formula} that we use for convergence analyses.

The structure of the primal objective function Hessian $\bbH_{k}$,
the degree matrix $\bbD$, and the oriented and unoriented
Laplacians $\bbL_o$ and $\bbL_u$ make distributed implementation
of \eqref{x_update_formula} possible. Indeed, the matrix
$2c\bbD+\bbH_k$ is block diagonal and its $i$-th diagonal block is
given by $2cd_{i}\bbI+\nabla^2f_{i}(\bbx_{i})$ which is locally
available for node $i$. Likewise, the inverse matrix
$(2c\bbD+\bbH_k)^{-1}$ is block diagonal and locally computable
since the $i$-th diagonal block is
$(2cd_{i}\bbI+\nabla^2f_{i}(\bbx_{i}))^{-1}$. Computations of the
products $\bbL_u\bbx_{k}$ and $\bbL_o\bbx_{k+1}$ can be
implemented in a decentralized manner as well, since the Laplacian
matrices $\bbL_u$ and $\bbL_o$ are block neighbor sparse in the
sense that the ${(i,j)}$-th block is not null if and only if nodes
$i$ and $j$ are neighbors or $j=i$. Therefore, nodes can compute
their local parts for the products $\bbL_u\bbx_{k}$ and
$\bbL_o\bbx_{k+1}$ by exchanging information with their neighbors.
By defining components of the vector $\bbphi_k$ as
$\bbphi_k:=[\bbphi_{1,k},\dots,\bbphi_{n,k}]$, the update formula
in \eqref{x_update_formula} for the individual agents can then be
written block-wise as
 \begin{align}\label{x_local_update_formula}
\bbx_{i,k+1}=\ & \left(2cd_{i}\bbI+\nabla^2f_{i}(\bbx_{i,k})\right)^{-1}\Big[ cd_{i}\bbx_{i,k}+c\sum_{j\in \ccalN_i}\bbx_{j,k}
\nonumber \\
&\qquad +\nabla^2f_{i}(\bbx_{i,k})\bbx_{i,k}-\nabla f_i(\bbx_{i,k})-\bbphi_{i,k} \Big],
\end{align}
where $\bbx_{i,k}$ corresponds to the iterate of node $i$ at step $k$. Notice that the defintion $\bbL_u:=(1/2)\bbE_u^T\bbE_u=(1/2)(\bbA_s+\bbA_d)^T(\bbA_s+\bbA_d)$ is used to simplify the $i$-th component of $c\bbL_u\bbx_{k}$ as $c\sum_{j\in \ccalN_i}(\bbx_{i,k}+\bbx_{j,k})$ which is equivalent to $cd_{i}\bbx_{i,k}+c\sum_{j\in \ccalN_i}\bbx_{j,k}
$. Further, using the definition $\bbL_o=(1/2)\bbE_o^T\bbE_o=(1/2)(\bbA_s-\bbA_d)^T(\bbA_s-\bbA_d)$, the $i$-th component of the product $c\bbL_o\bbx_{k+1}$ in \eqref{phi_local_update_formula} can be simplified as $c\sum_{j\in \ccalN_i}(\bbx_{i,k}-\bbx_{j,k})$. Therefore, the second update formula in \eqref{x_update_formula} can be locally implemented at each node $i$ as
\begin{equation}\label{phi_local_update_formula}
\bbphi_{i,k+1} =\bbphi_{i,k} +c \sum_{j\in \ccalN_i} \left(\bbx_{i,k+1} -\bbx_{j,k+1}\right).
\end{equation}
The proposed DQM method is summarized in Algorithm \ref{algo_DQM}. The initial value for the local iterate $\bbx_{i,0}$ can be any arbitrary vector in $\reals^p$. The initial vector $\bbphi_{i,0}$ should be in column space of $\bbE_o^T$. To guarantee satisfaction of this condition, the initial vector is set as $\bbphi_{i,0}=\bb0$. At each iteration $k$, updates of the primal and dual variables in \eqref{x_local_update_formula} and \eqref{phi_local_update_formula} are computed in Steps 2 and 4, respectively. Nodes exchange their local variables $\bbx_{i,k}$ with their neighbors $j\in \ccalN_i$ in Step 3, since this information is required for the updates in Steps 2 and 4.

%
\begin{algorithm}[t]{\small
\caption{DQM method at node $i$}\label{algo_DQM}
\begin{algorithmic}[1] {
\REQUIRE  Initial local iterates $\bbx_{i,0}$ and $\bbphi_{0}$.
\FOR {$k=0,1,2,\ldots$}

   \STATE Update the local iterate $\bbx_{i,k+1}$ as    \vspace{-2mm}
    \begin{align}
\bbx_{i,k+1}=\ & \left(2cd_{i}\bbI+\nabla^2f_{i}(\bbx_{i,k})\right)^{-1}\bigg[ cd_{i}\bbx_{i,k}+c\sum_{j\in \ccalN_i}\bbx_{j,k}
\nonumber \\
&\qquad\qquad+\nabla^2f_{i}(\bbx_{i,k})\bbx_{i,k}-\nabla f_i(\bbx_{i,k})-\bbphi_{i,k} \bigg].\nonumber
\end{align}

   \STATE Exchange iterates $\bbx_{i,k+1}$ with neighbors $\displaystyle{j\in \mathcal{N}_i}$.
   \STATE Update local dual variable $\bbphi_{k+1} $ as \\
         $\displaystyle{
\qquad \bbphi_{i,k+1} =\bbphi_{i,k} +c \sum_{j\in \ccalN_i} \left(\bbx_{i,k+1} -\bbx_{j,k+1}\right).}
$
\ENDFOR}
\end{algorithmic}}\end{algorithm}

%
DADMM, DQM, and DLM occupy different points in a tradeoff curve of
computational cost per iteration and number of iterations needed
to achieve convergence. The computational cost of each DADMM
iteration is large in general because it requires solution of the
optimization problem in \eqref{ADMM_x_update}. The cost of DLM
iterations is minimal because the solution of
\eqref{DLM_optimality_cond} can be reduced to the inversion of a
block diagonal matrix; see \cite{ling2014dlm}. The cost of DQM
iterations is larger than the cost of DLM iterations because they
require evaluation of local Hessians as well as inversion of the
matrices $2cd_{i}\bbI+\nabla^2f_{i}(\bbx_{i,k})$ to implement
\eqref{x_local_update_formula}. But the cost is smaller than the
cost of DADMM iterations except in cases in which solving
\eqref{ADMM_x_update} is easy. In terms of the number of
iterations required until convergence, DADMM requires the least
and DLM the most. The foremost technical conclusions of the
convergence analysis presented in the following section are: (i)
convergence of DQM is strictly faster than convergence of DLM;
(ii) asymptotically in the number of iterations, the per iteration
improvements of DADMM and DQM are identical. It follows from these
observations that DQM achieves target optimality in a number of
iterations similar to DADMM but with iterations that are
computationally cheaper.

%% file: Con_Anlysis.tex

\section{Convergence Analysis} \label{sec:Analysis}

In this section we show that the sequence of iterates $\bbx_k$
generated by DQM converges linearly to the optimal argument
$\bbx^*=[\tbx^*;\dots;\tbx^*]$. As a byproduct of this analysis we
also obtain a comparison between the linear convergence
constants of DLM, DQM, and DADMM. To derive these results we
make the following assumptions.

\begin{assumption}\label{network_eigenvalues}
The network is such that any singular value of the unoriented
incidence matrix $\bbE_u$, defined as $\sigma(\bbE_u)$, satisfies
$0<\gamma_u\leq\sigma(\bbE_u)\leq\Gamma_u$ where $\gamma_u$ and
$\Gamma_u$ are constants; the smallest non-zero singular value of
the oriented incidence matrix $\bbE_o$ is $\gamma_o>0$.
\end{assumption}

\begin{assumption}\label{convexity_assumption}
The local objective functions $f_i(\bbx)$ are twice differentiable
and the eigenvalues of their local Hessians $ \nabla^2 f_i(\bbx)$
are bounded within positive constants $m$ and $M$ where $0<m\leq
M<\infty$ so that for all $\bbx\in \reals^p$ it holds
\begin{equation}\label{local_hessian_eigenvlaue_bounds}
m\bbI\ \preceq\ \nabla^2 f_i(\bbx)\ \preceq\ M\bbI.
\end{equation}
\end{assumption}


\begin{assumption}\label{Lipschitz_assumption} The local Hessians $\nabla^2 f_i(\bbx)$ are Lipschitz continuous with constant $L$ so that for all $\bbx, \hbx \in \reals^p$ it holds
\begin{equation}
   \left\| \nabla^2 f_i(\bbx)-\nabla^2 f_i(\hbx) \right\| \ \leq\  L\ \| \bbx- \hbx \|.
\end{equation}

\end{assumption}


The eigenvalue bounds in Assumption \ref{network_eigenvalues} are measures of network
connectivity. Note that the assumption that all the singular values of the unoriented incidence matrix $\bbE_u$ are positive implies that the graph is non-bipartite. The conditions imposed by assumptions
\ref{convexity_assumption} and \ref{Lipschitz_assumption} are
typical in the analysis of second order methods; see, e.g.,
\cite[Chapter 9]{boyd}. The lower bound for the eigenvalues of the
local Hessians $\nabla^2 f_i(\bbx)$ implies strong convexity of
the local objective functions $f_{i}(\bbx)$ with constant $m$,
while the upper bound $M$ for the eigenvalues of the local
Hessians $\nabla^2 f_i(\bbx)$ is tantamount to Lipschitz
continuity of local gradients $\nabla f_i(\bbx)$ with Lipschitz
constant $M$. Further note that as per the definition of the
aggregate objective $f(\bbx):= \sum_{i=1}^{n}f_i(\bbx_i)$, the
Hessian $\bbH(\bbx):=\nabla^2f(\bbx)\in \reals^{np\times np}$ is
block diagonal with $i$-th diagonal block given by the $i$-th
local objective function Hessian $\nabla^2f_i(\bbx_i)$. Therefore,
the bounds for the local Hessians' eigenvalues in
\eqref{local_hessian_eigenvlaue_bounds} also hold for the
aggregate function Hessian. Thus, we have that for any $\bbx\in
\reals^{np}$ the eigenvalues of the Hessian $\bbH(\bbx)$ are
uniformly bounded as
\begin{equation}\label{aggregate_hessian_eigenvlaue_bounds}
m\bbI\ \preceq\ \bbH(\bbx)\ \preceq\ M\bbI.
\end{equation}
Assumption \ref{Lipschitz_assumption} also implies an analogous condition for the aggregate function Hessian $\bbH(\bbx)$ as we show in the following lemma.

%
\begin{lemma}\label{Hessian_Lipschitz_countinous}
Consider the definition of the aggregate function $f(\bbx):=
\sum_{i=1}^{n}f_i(\bbx_i)$. If Assumption
\ref{Lipschitz_assumption} holds true, the aggregate function
Hessian $\bbH(\bbx) =: \nabla^2f(\bbx)$ is Lipschitz continuous
with constant $L$. I.e., for all $\bbx,\hbx\in \reals^{np}$ we can
write
\begin{equation}\label{H_Lipschitz_claim}
\left\|\bbH(\bbx)-\bbH(\hbx)\right\| \leq L \| \bbx-\hbx\| .
\end{equation}
\end{lemma}

%
\begin{myproof}
See Appendix \ref{app_Hessian_Lipschitz}.
\end{myproof}



DQM can be interpreted as an attempt to approximate the primal update of DADMM. Therefore, we evaluate the performance of DQM by studying a measure of the error of the approximation in the DQM update relative to the DADMM update. In the primal update of DQM, the gradient $\nabla f(\bbx_{k+1})$ is estimated by the approximation $\nabla f(\bbx_{k})+ \bbH_k(\bbx_{k+1}-\bbx_{k})$. Therefore, we can define the DQM error vector $\bbe_{k}^{DQM}$ as
\begin{equation}\label{DQM_error}
\bbe_{k}^{DQM}:=\nabla f(\bbx_{k})+ \bbH_k(\bbx_{k+1}-\bbx_{k})-\nabla f(\bbx_{k+1}).
\end{equation}
Based on the definition in \eqref{DQM_error}, the approximation error of DQM vanishes when the difference of two consecutive iterates $\bbx_{k+1}-\bbx_{k}$ approaches zero. This observation is formalized in the following proposition by introducing an upper bound for the error vector norm $\|\bbe_k^{DQM}\|$ in terms of the difference norm $\|\bbx_{k+1}-\bbx_{k}\|$.

%
\begin{proposition}\label{error_vector_proposition}
Consider the DQM method as introduced in \eqref{DQM_update_1}-\eqref{DQM_update_2} and the error $\bbe_k^{DQM}$ defined in \eqref{DQM_error}. If Assumptions \ref{initial_val_assum}-\ref{Lipschitz_assumption} hold true, the DQM error norm $\| \bbe_k^{DQM}\|$ is bounded above by
\begin{equation}\label{bound_for_DQM_error}
    \left\|\bbe_{k}^{DQM} \right\| \leq  \min\left\{2M\|\bbx_{k+1}-\bbx_{k}  \|,\frac{L}{2}\|\bbx_{k+1}-\bbx_{k} \|^2\right\}.
\end{equation}
\end{proposition}

%
\begin{myproof} See Appendix  \ref{app_error_bound}. \end{myproof}


Proposition \ref{error_vector_proposition} asserts that the error
norm $\|\bbe_{k}^{DQM} \| $ is bounded above by the minimum of a
linear and a quadratic term of the iterate difference norm
$\|\bbx_{k+1}-\bbx_{k}  \|$. Hence, the approximation error
vanishes as the sequence of iterates $\bbx_k$ converges. We will
show in Theorem \ref{DQM_convergence} that the sequence
$\|\bbx_{k+1}-\bbx_k\|$ converges to zero which implies that the
error vector $\bbe_k^{DQM}$ converges to the null vector $\bb0$.
Notice that after a number of iterations the term
$(L/2)\|\bbx_{k+1}-\bbx_{k}\|$ becomes smaller than $2M$, which
implies that the upper bound in \eqref{bound_for_DQM_error} can be
simplified as $(L/2)\|\bbx_{k+1}-\bbx_{k}  \|^2$ for sufficiently
large $k$. This is important because it implies that the error
vector norm $ \|\bbe_{k}^{DQM} \|$ eventually becomes proportional
to the quadratic term $\|\bbx_{k+1}-\bbx_k\|^2$ and, as a
consequence, it vanishes faster than the term
$\|\bbx_{k+1}-\bbx_k\|$.

Utilize now the definition in \eqref{DQM_error} to rewrite the
primal variable DQM update in \eqref{DQM_update_1} as 
\begin{align}
\nabla f(\bbx_{k+1})+\bbe_{k}^{DQM}+\bbA^T\bblambda_{k}+c \bbA^{T}\! \left( \bbA\bbx_{k+1}\!+\!\bbB \bbz_{k}  \right)&=\bb0.\label{joint_1}
\end{align}
Comparison of \eqref{joint_1} with the optimality condition for
the DADMM update in \eqref{ADMM_x_update} shows that they coincide
except for the  gradient approximation error term
$\bbe_{k}^{DQM}$. The DQM and DADMM updates for the auxiliary
variables $\bbz_k$ and the dual variables $\bblambda_k$ are
identical [cf. \eqref{ADMM_z_update}, \eqref{ADMM_lambda_update},
and \eqref{DQM_update_2}], as already observed.

Further let the pair $(\bbx^*,\bbz^*)$ stand for the unique
solution of \eqref{original_optimization_problem2} with uniqueness
implied by the strong convexity assumption and define $\bbalpha^*$
as the unique optimal multiplier that lies in the column space of
$\bbE_o$ -- see Lemma 1 of \cite{cQingRibeiroADMM14} for a proof
that such optimal dual variable exists and is unique. To study
convergence properties of DQM we modify the system of DQM
equations defined by \eqref{DQM_update_2} and \eqref{joint_1},
which is equivalent to the system \eqref{DQM_update_1} --
\eqref{DQM_update_2}, to include terms that involve differences
between current iterates and the optimal arguments $\bbx^*$,
$\bbz^*$, and $\bbalpha^*$. We state this reformulation in the
following lemma.

\begin{lemma}\label{equalities_for_optima_lemma}
Consider the DQM method as defined by \eqref{DQM_update_1}-\eqref{DQM_update_2} and its equivalent formulation in \eqref{DQM_update_2} and \eqref{joint_1}. If Assumption \ref{initial_val_assum} holds true, then the optimal arguments $\bbx^*$, $\bbz^*$, and $\bbalpha^*$ satisfy
\begin{align}
   \nabla f(\bbx_{k+1})-\nabla f(\bbx^*)+\bbe_{k}^{DQM}
          + \bbE_o^T(\bbalpha_{k+1}-\bbalpha^*) \nonumber \\
         -c\bbE_u^T \left( \bbz_{k}-\bbz_{k+1}  \right)
   & = \bb0, \label{relation1} \\
   2(\bbalpha_{k+1}-\bbalpha_k)-{c}\bbE_o(\bbx_{k+1}-\bbx^*)
   &=\bb0, \label{relation2} \\
   \bbE_u(\bbx_{k}-\bbx^*)-2(\bbz_k-\bbz^*)
   &=\bb0.\label{relation3}
\end{align}
\end{lemma}
\begin{myproof}
See Appendix \ref{app_equalities_for_optima}.
\end{myproof}


With the preliminary results in Lemmata \ref{Hessian_Lipschitz_countinous} and \ref{equalities_for_optima_lemma} and Proposition \ref{error_vector_proposition} we can state our convergence results. To do so, define the energy function $V: \reals^{mp\times mp}\to\reals$ as
\begin{equation}\label{eqn_energy_function}
V(\bbz,\bbalpha):=c\|\bbz-\bbz^{*}\|^2+\frac{1}{c}\|\bbalpha-\bbalpha^*\|^2.
\end{equation}
The energy function $V(\bbz,\bbalpha)$ captures the distances of the variables $\bbz_{k}$ and $\bbalpha_{k}$ to the respective optimal arguments $\bbz^{*} $ and $\bbalpha^*$. To simplify notation we further define the variable $\bbu\in \reals^{2mp}$ and matrix $\bbC\in \reals^{2mp\times 2mp}$ as
\begin{equation}\label{C_u_definitions}
\bbu:= \left[ \begin{array}{rr}
\bbz  \\
\bbalpha  \end{array} \right], \quad
\bbC:= \left[ \begin{array}{rr}
c\bbI_{mp} & \bb0  \\
\bb0 & (1/c)\bbI_{mp} \end{array} \right].
\end{equation}
Based on the definitions in \eqref{C_u_definitions}, the energy function in \eqref{eqn_energy_function} can be alternatively written $V(\bbz,\bbalpha) = V(\bbu) = \|\bbu-\bbu^*\|_{\bbC}^2$, where $\bbu^*=[\bbz^*;\bbalpha^*]$. The energy sequence $V(\bbu_k)=\|\bbu_{k}-\bbu^*\|_\bbC^2$ converges to zero at a linear rate as we state in the following theorem.

\begin{thm}\label{DQM_convergence}
Consider the DQM method as defined by \eqref{DQM_update_1}-\eqref{DQM_update_2}, let the constant $c$ be such that $c> 4M^2/({m\gamma_u^2})$, and define the sequence of non-negative variables $\zeta_k$ as
\begin{equation}\label{zeta_DQM}
\zeta_k:=\min\left\{\frac{L}{2}\|\bbx_{k+1}-\bbx_{k} \|, 2M\right\}.
\end{equation}
Further, consider arbitrary constants $\mu$, $\mu'$, and $\eta$ with $\mu,\mu'>1$ and $\eta_k\in({\zeta_k}/{m},{c\gamma_u^2}/{\zeta_k})$. If Assumptions \ref{initial_val_assum}-\ref{Lipschitz_assumption} hold true, then the sequence $\|\bbu_{k}-\bbu^*\|_{\bbC}^2$ generated by DQM satisfies
\begin{equation}\label{DQM_linear_claim}
\|\bbu_{k+1}-\bbu^*\|_{\bbC}^2\  \leq\ \frac{1}{1+\delta_k}  \|\bbu_{k}-\bbu^*\|_{\bbC}^2\ ,
\end{equation}
where the sequence of positive scalars $\delta_k$ is given by
\begin{equation}\label{DLM_DQM_delta}
\delta_k=\min
\Bigg\{
\frac{(\mu-1)(c\gamma_u^2-\eta_k\zeta_k)\gamma_o^{2}}
     {\mu\mu'( c\Gamma_u^2\gamma_u^2+4\zeta_k^2/c(\mu'-1))}
,
\frac{  m-{\zeta_k}/{\eta_k} }{{ c}\Gamma_u^2/4
+\mu M^2/c\gamma_o^{2}}
\Bigg\}.
\end{equation}
\end{thm}

\begin{myproof}
See Appendix \ref{linear_convg_app}.
\end{myproof}


Notice that $\delta_k$ is a decreasing function of $\zeta_k$ and that $\zeta_k$ is bounded above by $2M$. Therefore, if we substitute $\zeta_k$ by $2M$ in \eqref{DLM_DQM_delta}, the inequality in \eqref{DQM_linear_claim} is still valid. This substitution implies that the sequence $\|\bbu_{k}-\bbu^*\|_{\bbC}^2$ converges linearly to zero with a coefficient not larger than $1-\delta$ with $\delta=\delta_k$ following from \eqref{DQM_linear_claim} with $\zeta_k=2M$. The more generic definition of $\zeta_k$ in \eqref{zeta_DQM} is important for the rate comparisons in Section \ref{sec:rate_comparison}. Observe that in order to guarantee that $\delta_k>0$ for all $k\geq 0 $, $\eta_k$ is chosen from the interval $({\zeta_k}/{m},{c\gamma_u^2}/{\zeta_k})$. This interval is non-empty since the constant $c$ is chosen as $c> {4M^2}/({m\gamma_u^2})\geq  {\zeta_k^2}/({m\gamma_u^2})$.

The linear convergence in Theorem \ref{DQM_convergence} is for the vector $\bbu_k$ which includes the auxiliary variable $\bbz_k$ and the multipliers $\bbalpha_k$. Linear convergence of the primal variables $\bbx_k$ to the optimal argument $\bbx^*$ follows as a corollary that we establish next.

\begin{corollary}\label{convergence_of_primal}
Under the assumptions in Theorem \ref{DQM_convergence}, the sequence of squared norms $\|\bbx_k-\bbx^*\|^2$ generated by the DQM algorithm converges R-linearly to zero, i.e.,
\begin{equation}\label{r_linear_cliam}
\|\bbx_k-\bbx^*\|^2\leq \frac{4}{c\gamma_u^2}\|\bbu_k-\bbu^*\|_\bbC^2.
\end{equation}
\end{corollary}

\begin{myproof}
Notice that according to \eqref{relation3} we can write $\|\bbE_u(\bbx_k-\bbx^*)\|^2=4\|\bbz_k-\bbz^*\|^2$. Since $\gamma_u$ is the smallest singular value of $\bbE_u$, we obtain that $\|\bbx_k-\bbx^*\|^2\leq({4}/{\gamma_u^2})\|\bbz_k-\bbz^*\|^2.$ Moreover, according to the relation $\|\bbu_k-\bbu^*\|_\bbC^2=c\|\bbz_k-\bbz^{*}\|^2+({1}/{c})\|\bbalpha_k-\bbalpha^*\|^2$ we can write  $c\|\bbz_k-\bbz^*\|^2\leq\|\bbu_k-\bbu^*\|_\bbC^2.$ Combining these two inequalities yields the claim in \eqref{r_linear_cliam}.
\end{myproof}


As per Corollary \ref{convergence_of_primal}, convergence of the sequence $\bbx_k$ to $\bbx^*$ is dominated by a linearly decreasing sequence. Notice that the sequence of squared norms $\|\bbx_k-\bbx^*\|^2$ need not be monotonically decreasing as the energy sequence $\|\bbu_{k+1}-\bbu^*\|_{\bbC}^2$ is.

\subsection{Convergence rates comparison}\label{sec:rate_comparison}

Based on the result in Corollary \ref{convergence_of_primal}, the sequence of iterates $\bbx_k$ generated by DQM converges. This observation implies that the sequence $\|\bbx_{k+1}-\bbx_k\|$ approaches zero. Hence, the sequence of scalars $\zeta_k$ defined in \eqref{zeta_DQM} converges to 0 as time passes, since $\zeta_k$ is bounded above by $(L/2)\|\bbx_{k+1}-\bbx_k\|$. Using this fact that $\lim_{k\to \infty }\zeta_k=0$ to compute the limit of $\delta_k$ in \eqref{DLM_DQM_delta} and further making $\mu'\rightarrow1$ in the resulting limit we have that
\begin{align}\label{DQM_delta}
\lim_{k\to \infty} \delta_k=\min
&\Bigg\{
\frac{ (\mu-1)\gamma_o^2}{ {\mu\Gamma_u^2}}
,\
\frac{  m}{c\Gamma_u^2/4
+\mu M^2/c\gamma_o^{2}}
\Bigg\}.
\end{align}
Notice that the limit of $\delta_k$ in \eqref{DQM_delta} is
identical to the constant of linear convergence for DADMM
\cite{Shi2014-ADMM}. Therefore, we conclude that as time passes
the constant of linear convergence for DQM approaches the one
for DADMM.

To compare the convergence rates of DLM, DQM and DADMM we define the error of the gradient approximation for DLM as
\begin{equation}\label{DLM_error}
\bbe_{k}^{DLM}=\nabla f(\bbx_{k})+ \rho(\bbx_{k+1}-\bbx_{k})-\nabla f(\bbx_{k+1}),
\end{equation}
which is the difference of exact gradient $\nabla f(\bbx_{k+1})$ and the DLM gradient approximation $\nabla f(\bbx_{k})+ \rho(\bbx_{k+1}-\bbx_{k})$. Similar to the result in Proposition \ref{error_vector_proposition} for DQM we can show that the DLM error vector norm $\| \bbe_k^{DLM}\|$ is bounded by a factor of $\|\bbx_{k+1}-\bbx_{k}  \|$.
%
%
\begin{proposition}\label{error_vector_proposition_2}
Consider the DLM algorithm with updates in \eqref{ADMM_z_update}-\eqref{DLM_x_update} and the error vector $\bbe_k^{DLM}$ defined in \eqref{DLM_error}. If Assumptions \ref{initial_val_assum}-\ref{Lipschitz_assumption} hold true,  the DLM error vector norm $\| \bbe_k^{DLM}\|$ satisfies
\begin{equation}\label{bound_for_DLM_error}
    \left\|\bbe_{k}^{DLM} \right\| \leq (\rho+M)\|\bbx_{k+1}-\bbx_{k}  \|.
\end{equation}
\end{proposition}

%
\begin{myproof} See Appendix  \ref{app_error_bound}. \end{myproof}

The result in Proposition \ref{error_vector_proposition_2} differs from Proposition \ref{error_vector_proposition} in that the DLM error $\|\bbe_{k}^{DLM}\|$ vanishes at a rate of $\|\bbx_{k+1}-\bbx_{k}\|$ whereas the DQM error $\|\bbe_{k}^{DQM}\|$ 
eventually becomes proportional to $\|\bbx_{k+1}-\bbx_k\|^2$. This results in DLM failing to approach the convergence behavior of DADMM as we show in the following theorem.

\begin{thm}\label{DLM_convergence}
Consider the DLM method as introduced in \eqref{ADMM_z_update}-\eqref{DLM_x_update}.
Assume that the constant $c$ is chosen such that $c> (\rho+M)^2/({m\gamma_u^2})$.
Moreover, consider $\mu,\mu'>1$ as arbitrary constants and $\eta$ as a positive constant chosen from the interval $((\rho+M)/{m},{c\gamma_u^2}/{(\rho+M)})$.
If Assumptions \ref{initial_val_assum}-\ref{Lipschitz_assumption} hold true, then the sequence $\|\bbu_{k}-\bbu^*\|_{\bbC}^2$ generated by DLM satisfies
\begin{equation}\label{DLM_linear_claim}
\|\bbu_{k+1}-\bbu^*\|_{\bbC}^2\  \leq\ \frac{1}{1+\delta}  \|\bbu_{k}-\bbu^*\|_{\bbC}^2\ ,
\end{equation}
where the scalar $\delta$ is given by
\begin{equation}\label{DLM_delta}
\delta\!=\!\min\!
\Bigg\{\!
\frac{(\mu-1)(c\gamma_u^2-\eta_k(\rho\!+\!M))\gamma_o^{2}}
     {\mu\mu'( c\Gamma_u^2\gamma_u^2\!+\!4(\rho\!+\!M)^2\!/c(\mu'\!-\!1))}
,
\frac{  m-{(\rho\!+\!M)}/{\eta_k} }{{ c}\Gamma_u^2/4
\!+\!\mu M^2\!/c\gamma_o^{2}}
\!\Bigg\}
\end{equation}
\end{thm}
\begin{myproof}
See Appendix \ref{linear_convg_app}.
\end{myproof}


Based on the result in Theorem \ref{DLM_convergence}, the sequence
$\|\bbu_{k+1}-\bbu^*\|_{\bbC}^2$ generated by DLM converges
linearly to 0. This result is similar to the convergence
properties of DQM as shown in Theorem \ref{DQM_convergence};
however, the constant of linear convergence $1/(1+\delta)$ in
\eqref{DLM_linear_claim} is smaller than the constant
$1/(1+\delta_k)$ in \eqref{DQM_delta}.

%% file: Simulations.tex

%
\section{Numerical analysis}\label{sec:simulations}

In this section we compare the performances of DLM, DQM and DADMM
in solving a logistic regression problem. Consider a training set
with points whose classes are known and the goal is finding the
classifier that minimizes the loss function. Let $q$ be the number
of training points available at each node of the network.
Therefore, the total number of training points is $nq$. The
training set $\{\bbs_{il}, y_{il}\}_{l=1}^q$ at node $i$ contains
$q$ pairs of $(\bbs_{il}, y_{il})$, where $\bbs_{il}$ is a feature
vector and $y_{il}\in \{-1,1\}$ is the corresponding class. The
goal is to estimate the probability $\Pc{y=1\mid \bbs}$ of having
label $y=1$ for a given feature vector $\bbs$ whose class is not
known. Logistic regression models this probability as $\Pc{y=1\mid
\bbs}=1/(1+\exp(-\bbs^T\tbx))$ for a linear classifier $\tbx$ that
is computed based on the training samples. It follows from this
model that the maximum log-likelihood estimate of the classifier
$\tbx$ given the training samples
$\{\{\bbs_{il},y_{il}\}_{l=1}^q\}_{i=1}^n$ is
\begin{align}\label{eqn_logistic_regrssion_max_likelihood}
   \tbx^* \ :=\  \argmin_{\tbx\in \reals^p }  \    \sum_{i=1}^n \sum_{l=1}^{q}
                                  \log \Big[1+\exp(-y_{il}\bbs_{il}^T\tbx)\Big].
\end{align}
%
The optimization problem in
\eqref{eqn_logistic_regrssion_max_likelihood} can be written in
the form \eqref{original_optimization_problem1}. To do so, simply
define the local objective functions $f_i$ as
\begin{equation}
   f_i(\tbx) =     \sum_{l=1}^{q} \log \Big[1+\exp(-y_{il}\bbs_{il}^T\tbx)\Big].
\end{equation}
%
We define the optimal argument for decentralized optimization as
$\bbx^*=[\tbx^*;\dots;\tbx^*]$. Note that the reference (ground
true) logistic classifiers $\tbx^*$ for all the experiments in
this section are pre-computed with a centralized method.

\subsection{Comparison of DLM, DQM, and DADMM}\label{comparison}

%
\begin{figure}[t]
\centering
\vspace{-2mm}
\includegraphics[width=\linewidth,height=0.55\linewidth]{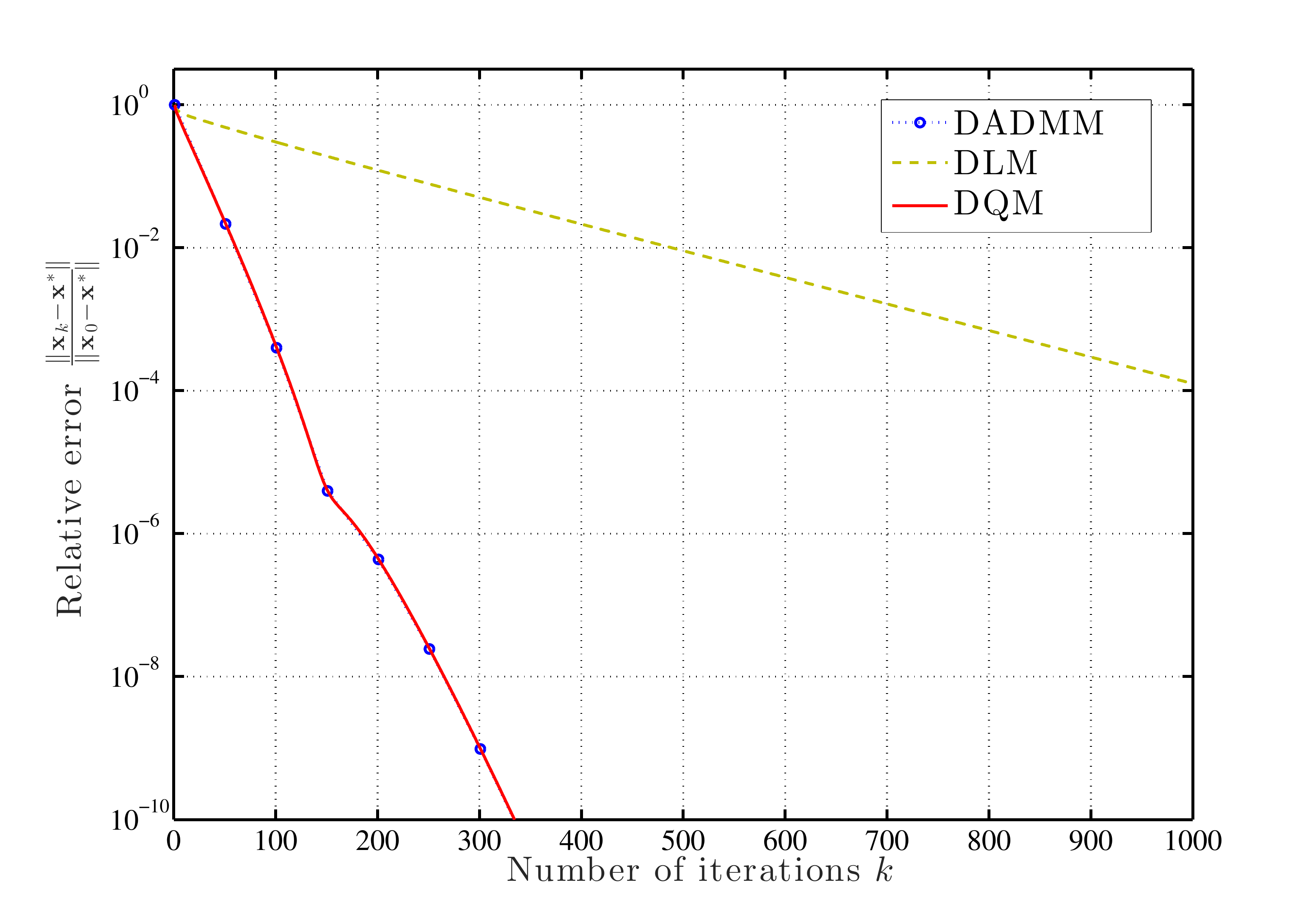}
\caption{{Relative error ${\|\bbx_k-\bbx^*\|}/{\|\bbx_0-\bbx^*\|}$ of DADMM, DQM, and DLM versus number of iterations for a random network of size $n=10$. The convergence path of DQM is similar to the one for DADMM and they  outperform DLM by orders of magnitude.
\vspace{-7mm}
}}
\label{fig2}
\end{figure}
%
\begin{figure}[t]
\centering
\vspace{-0.6mm}
\includegraphics[width=\linewidth,height=0.55\linewidth]{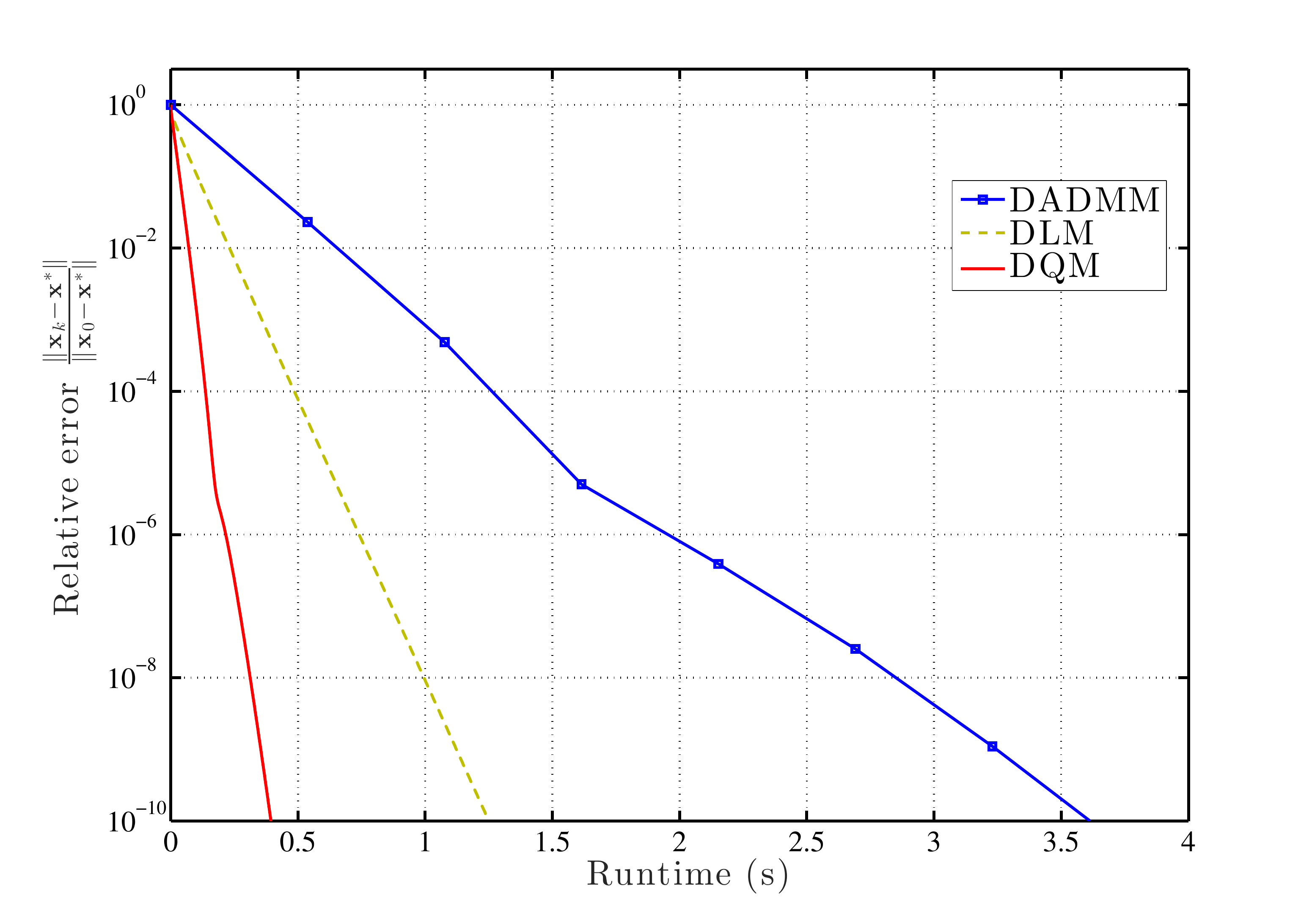}
\caption{{Relative error ${\|\bbx_k-\bbx^*\|}/{\|\bbx_0-\bbx^*\|}$
of DADMM, DQM, and DLM versus runtime for the setting in Fig.
\ref{fig2}. The computational cost of DQM is lower than DADMM and
DLM. }} \label{fig3}
\end{figure}

We compare the convergence paths of the DLM, DQM, and DADMM algorithms for solving the logistic regression problem in \eqref{eqn_logistic_regrssion_max_likelihood}. Edges between the nodes are randomly generated with the connectivity ratio $r_c$. Observe that the connectivity ratio $r_c$ is the probability of two nodes being connected.

In the first experiment we set the number of nodes as $n=10$ and
the connectivity ratio as $r_{\mathrm{c}}=0.4$. Each agent holds
$q=5$ samples and the dimension of feature vectors is $p=3$. Fig.
\ref{fig2} illustrates the relative errors
${\|\bbx_k-\bbx^*\|}/{\|\bbx_0-\bbx^*\|}$ for DLM, DQM, and DADMM
versus the number of iterations. Notice that the parameter $c$ for
the three methods is optimized by $c_\mathrm{ADMM}=0.7$,
$c_\mathrm{DLM}=5.5$, and $c_\mathrm{DQM}=0.7$. The convergence
path of DQM is almost identical to the convergence path of DADMM.
Moreover, DQM outperforms DLM by orders of magnitude. To be more
precise, the relative errors
${\|\bbx_k-\bbx^*\|}/{\|\bbx_0-\bbx^*\|}$ for DQM and DADMM after
$k=300$ iterations are below $10^{-9}$, while for DLM the relative
error after the same number of iterations is $5\times 10^{-2}$.
Conversely, achieving accuracy
${\|\bbx_k-\bbx^*\|}/{\|\bbx_0-\bbx^*\|}=10^{-3}$ for DQM and
DADMM requires $91$ iterations, while DLM requires $758$
iterations to reach the same accuracy. Hence, the number of
iterations that DLM requires to achieve a specific accuracy is $8$
times more than the one for DQM.

%
\begin{figure}[t]
\centering
\vspace{-2mm}
\includegraphics[width=\linewidth,height=0.55\linewidth]{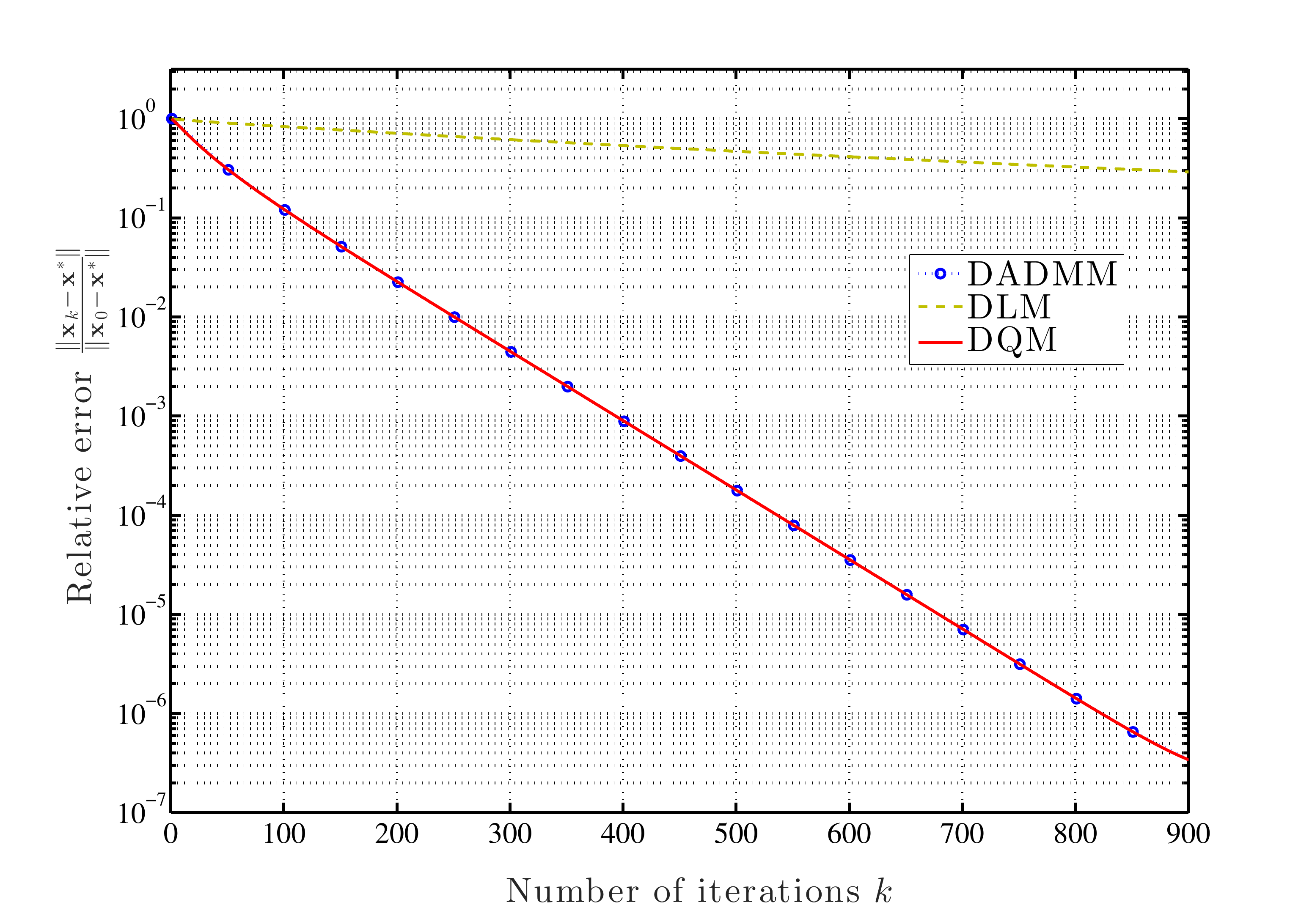}
\caption{Relative error ${\|\bbx_k-\bbx^*\|}/{\|\bbx_0-\bbx^*\|}$ of DADMM, DQM, and DLM versus number of iterations for a random network of size $n=100$. The performances of DQM and DADMM are still similar. DLM is impractical in this setting.
}
\vspace{-4mm}
\label{fig222}
\end{figure}

%
\begin{figure}[t]
\centering
\includegraphics[width=\linewidth,height=0.55\linewidth]{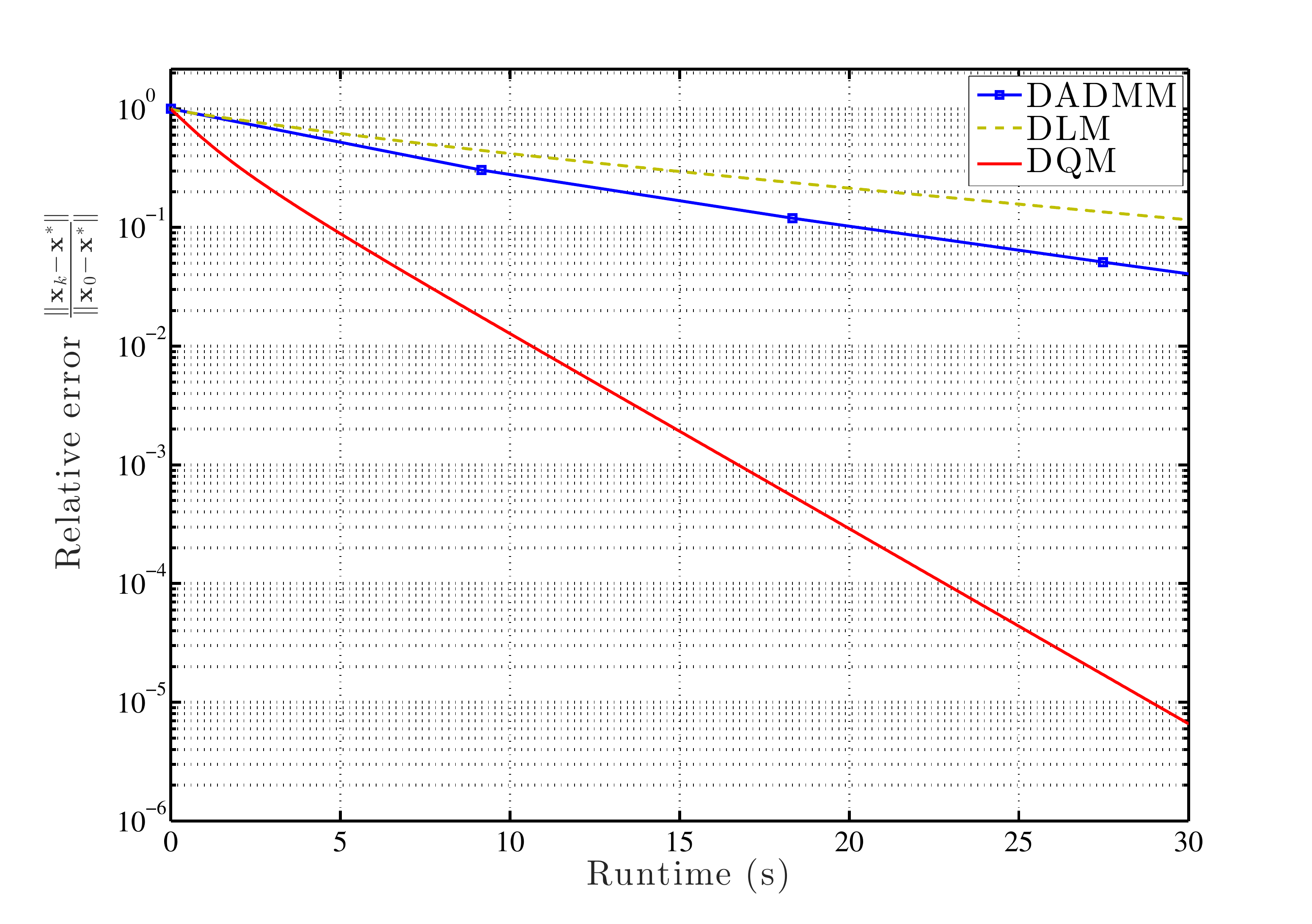}
\caption{{Relative error ${\|\bbx_k-\bbx^*\|}/{\|\bbx_0-\bbx^*\|}$ of DADMM, DQM, and DLM versus runtime for the setting in Fig. \ref{fig3}. The convergence time of DADMM is slightly faster relative to DLM, while DQM is the most efficient method among these three algorithms.
}}
\label{fig333}
\end{figure}

Observe that the computational complexity of DQM is lower than DADMM. Therefore, DQM outperforms DADMM in terms of convergence time or number of required operations until convergence. This phenomenon is shown in Fig \ref{fig3} by comparing the relative of errors of DLM, DQM, and DADMM versus CPU runtime. According to Fig \ref{fig3}, DADMM achieves the relative error ${\|\bbx_k-\bbx^*\|}/{\|\bbx_0-\bbx^*\|}=10^{-10}$ after running for $3.6$ seconds, while DLM and DQM require $1.3$ and $0.4$ seconds, respectively, to achieve the same accuracy.

We also compare the performances of DLM, DQM, and DADMM in a
larger scale logistic regression problem by setting size of
network $n=100$, number of sample points at each node $q=20$, and
dimension of feature vectors $p=10$. We keep the rest of the
parameters as in Fig. \ref{fig2}. Convergence paths of the
relative errors ${\|\bbx_k-\bbx^*\|}/{\|\bbx_0-\bbx^*\|}$ for DLM,
DQM, and DADMM versus the number of iterations are illustrated in
Fig. \ref{fig222}. Different choices of parameter $c$ are
considered for these algorithms and the best for each is chosen
for the final comparison. The optimal choices of parameter $c$ for
DADMM, DLM, and DQM are $c_\mathrm{ADMM}=0.68$,
$c_\mathrm{DLM}=12.3$, and $c_\mathrm{DQM}=0.68$, respectively.
The results for the large scale problem in Fig. \ref{fig222} are
similar to the results in Fig. \ref{fig2}. We observe that DQM
performs as well as DADMM, while both outperform DLM. To be more
precise, DQM and DADMM after $k=900$ iterations reach the relative
error ${\|\bbx_k-\bbx^*\|}/{\|\bbx_0-\bbx^*\|}=3.4\times10^{-7}$,
while the relative error of DLM after the same number of
iterations is $2.9\times 10^{-1}$. Conversely, achieving the
accuracy ${\|\bbx_k-\bbx^*\|}/{\|\bbx_0-\bbx^*\|}=0.3$ for DQM and
DADMM requires $52$ iterations, while DLM requires $870$
iterations to reach the same accuracy. Hence, in this setting the
number of iterations that DLM requires to achieve a specific
accuracy is $16$ times more than the one for DQM. These numbers
show that the advantages of DQM relative to DLM are more
significant in large scale problems.

Notice that in large scale logistic regression problems we expect
larger condition number for the objective function $f$. In these
scenarios we expect to observe a poor performance by the DLM
algorithm that only operates on first-order information. This
expectation is satisfied by comparing the relative errors of DLM,
DQM, and DADMM versus runtime for the large scale problem in Fig.
\ref{fig333}. In this case, DLM is even worse than DADMM that has
a very high computational complexity. Similar to the result in
Fig. \ref{fig222}, DQM has the best performance among these three
methods.

\subsection{Effect of the regularization parameter $c$}\label{dif_para}

\begin{figure}[t]
    \centering
    \vspace{-2mm}
\includegraphics[width=\linewidth,height=0.55\linewidth]{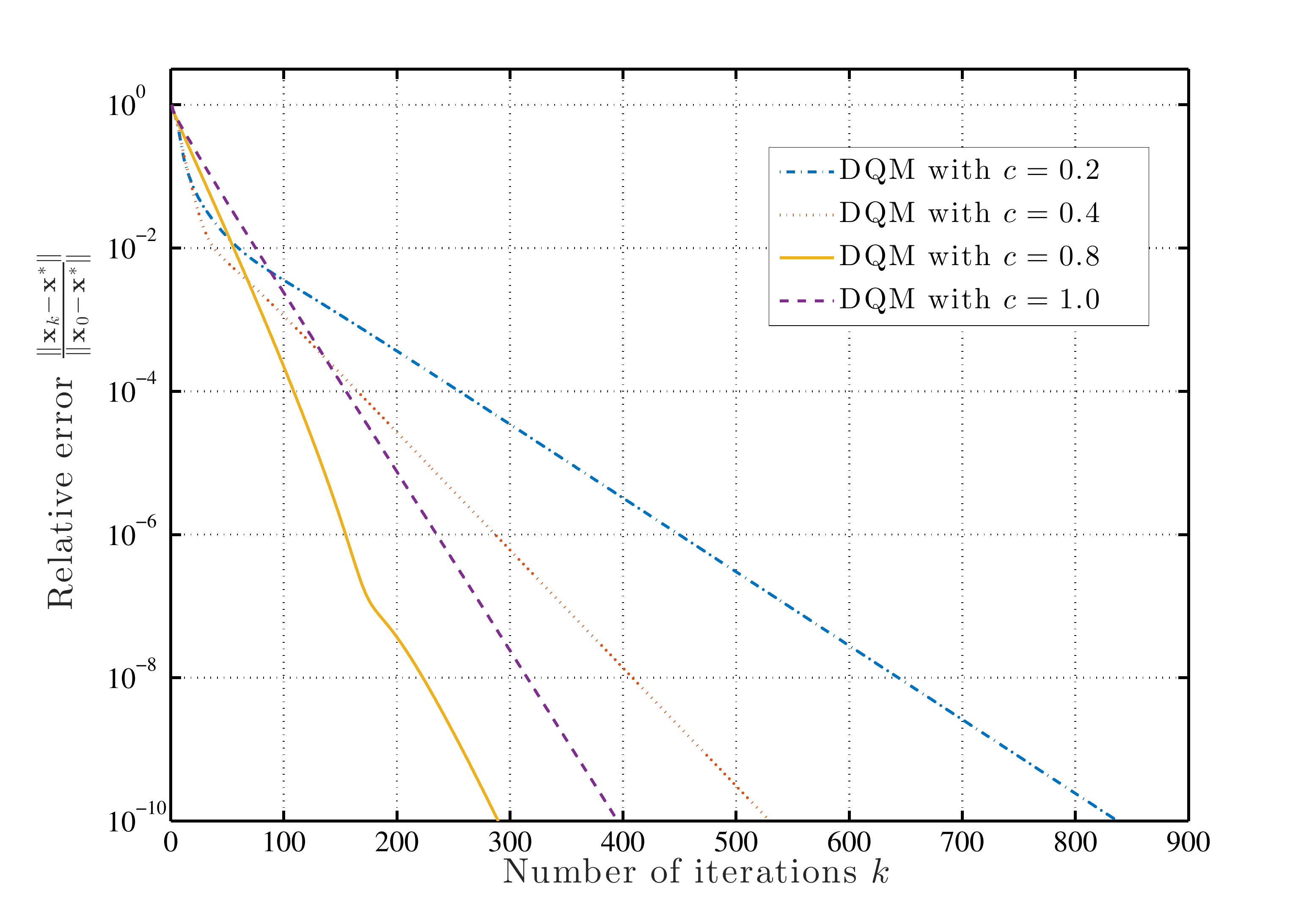}
\caption{Relative error ${\|\bbx_k-\bbx^*\|}/{\|\bbx_0-\bbx^*\|}$ of DQM for parameters $c=0.2$, $c=0.4$, $c=0.8$, and $c=1$ when the network is formed by $n=10$ nodes and the connectivity ratio is $r_c=0.4$. The best performance belongs to $c=0.8$.}\label{eps:Par_Tun}
\end{figure}

The parameter $c$ has a significant role in the convergence of
DADMM. Likewise, choosing the optimal choice of $c$ is critical in
the convergence of DQM. We study the effect of $c$ by tuning this
parameter for a fixed network and training set. We use all the
parameters in Fig. \ref{fig2} and we compare performance of the
DQM algorithm for the values $c=0.2$, $c=0.4$, $c=0.8$, and $c=1$.
Fig. \ref{eps:Par_Tun} illustrates the convergence paths of the
DQM algorithm for different choices of the parameter $c$. The best
performance among these choices is achieved for $c=0.8$. The
comparison of the plots in Fig. \ref{eps:Par_Tun} shows that
increasing or decreasing the parameter $c$ is not necessarily
leads to a faster convergence. We can interpret $c$ as the
stepsize of DQM which the optimal choice may vary for the problems
with different network sizes, network topologies, condition
numbers of objective functions, etc.

%

%
\begin{figure}[t]
    \centering
    \vspace{-2mm}
\includegraphics[width=\linewidth,height=0.55\linewidth]{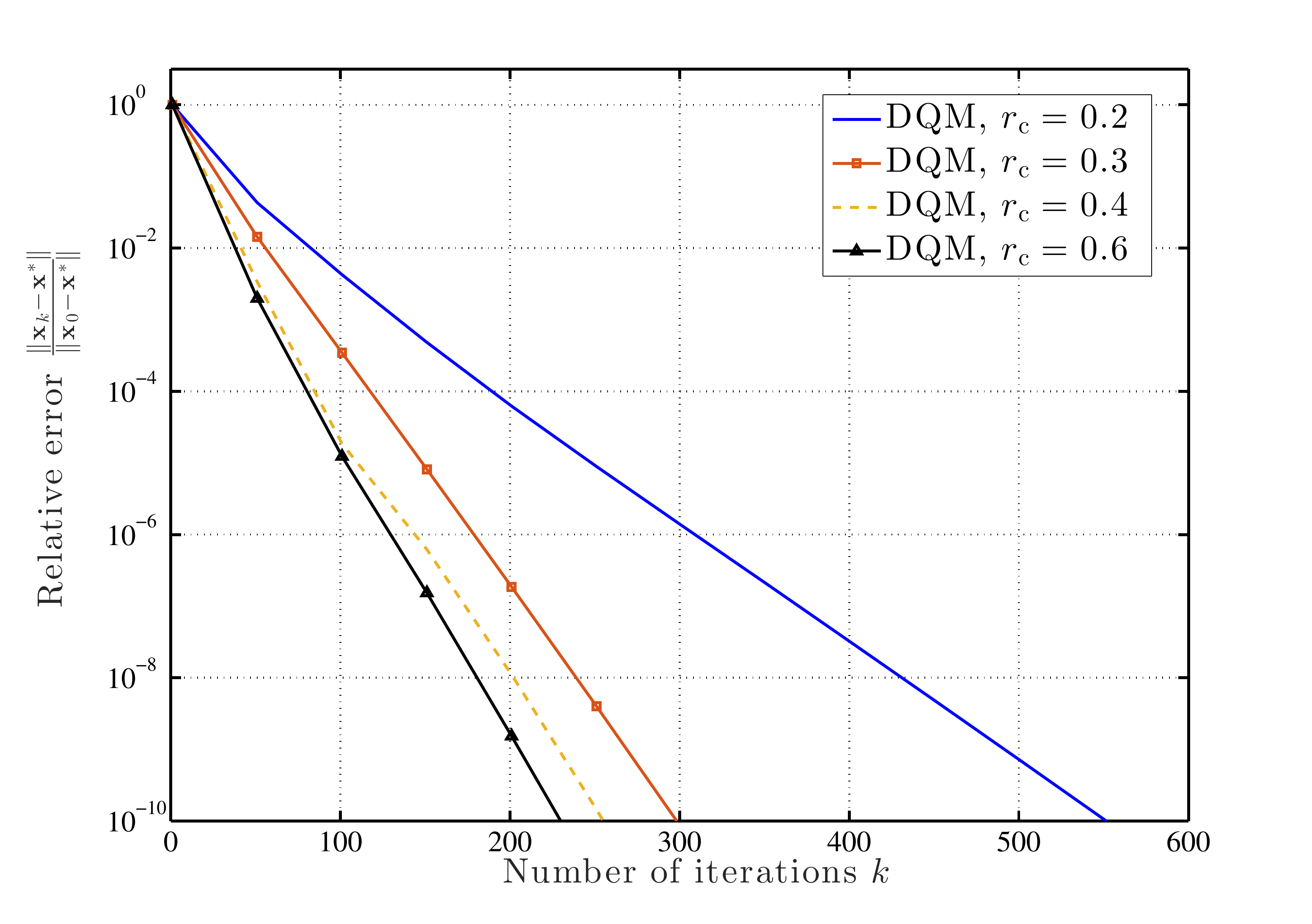}
  \caption{Relative error ${\|\bbx_k-\bbx^*\|}/{\|\bbx_0-\bbx^*\|}$ of DQM for random graphs with different connectivity ratios $r_c$. The linear convergence of DQM accelerates by increasing the connectivity ratio.}\label{eps:Net_Top}
\end{figure}

\subsection{Effect of network topology}\label{sec:Net_Top}

According to \eqref{DLM_DQM_delta} the constant of linear
convergence for DQM depends on the bounds for the singular values
of the oriented and unoriented incidence matrices $\bbE_o$ and
$\bbE_u$. These bounds are related to the connectivity ratio of
network. We study how the network topology affects the convergence
speed of DQM. We use different values for the connectivity ratio
to generate random graphs with different number of edges. In this
experiment we use the connectivity ratios
$r_c=\{0.2,0.3,0.4,0.6\}$ to generate the networks. The rest of
the parameters are the same as the parameters in Fig. \ref{fig2}.
Notice that since the connectivity parameters of these graphs are
different, the optimal choices of $c$ for these graphs are
different. The convergence paths of DQM with the connectivity
ratios $r_c=\{0.2,0.3,0.4,0.6\}$ are shown in Fig.
\ref{eps:Net_Top}. The optimal choices of the parameter $c$ for
these graphs are $c_\mathrm{0.2}=0.28$, $c_\mathrm{0.3}=0.25$,
$c_\mathrm{0.4}=0.31$, and $c_\mathrm{0.6}=0.28$, respectively.
Fig. \ref{eps:Net_Top} shows that the linear convergence of DQM
accelerates by increasing the connectivity ratio of the graph.

%

%% file: Conclusions.tex

\section{Conclusions}
A decentralized quadratically approximated version of the
alternating direction method of multipliers (DQM) is proposed for
solving decentralized optimization problems where components of
the objective function are available at different nodes of a
network. DQM minimizes a quadratic approximation of the convex
problem that DADMM solves exactly at each step, and hence reduces
the computational complexity of DADMM. Under some mild
assumptions, linear convergence of the sequence generated by DQM
is proven. Moreover, the constant of linear convergence for DQM
approaches that of DADMM asymptotically. Numerical results for a
logistic regression problem verify the analytical results that
convergence paths of DQM and DADMM are similar for large iteration
index, while the computational complexity of DQM is significantly
smaller than DADMM.

%% file: Appendix.tex

\begin{appendices}


\section{Proof of Lemma \ref{lag_var_lemma}\label{app_lag_var_lemma}}

According to the update for the Lagrange multiplier $\bblambda$ in \eqref{DQM_update_2}, we can substitute $\bblambda_k$ by $\bblambda_{k+1}-c\left( \bbA\bbx_{k+1}+\bbB\bbz_{k+1} \right)$. Applying this substitution into the first equation of \eqref{DQM_update_2} leads to
\begin{equation}\label{shajarian0}
\bbB^T\bblambda_{k+1}=\bb0.
\end{equation}
Observing the definitions  $\bbB=[-\bbI_{mp};-\bbI_{mp}]$ and $\bblambda=[\bbalpha;\bbbeta]$, and the result in \eqref{shajarian0}, we obtain $\bbalpha_{k+1}=-\bbbeta_{k+1}$ for $k\geq0$. Considering the initial condition $\bbalpha_{0}=-\bbbeta_{0}$, we obtain that $\bbalpha_{k}=-\bbbeta_{k}$ for $k\geq0$ which follows the first claim in Lemma \ref{lag_var_lemma}.

Based on the definitions $\bbA=[\bbA_s;\bbA_d]$, $\bbB=[-\bbI_{mp};-\bbI_{mp}]$, and $\bblambda=[\bbalpha;\bbbeta]$, we can split the update for the Lagrange multiplier $\bblambda$ in \eqref{ADMM_lambda_update} as
\begin{align}
\bbalpha_{k+1}&=\bbalpha_{k}+c[\bbA_s\bbx_{k+1}-\bbz_{k+1}] \label{sorry1},\\
\bbbeta_{k+1}&=\bbbeta_{k}+c[\bbA_d\bbx_{k+1}-\bbz_{k+1}] \label{sorry2}.
\end{align}
Observing the result that $\bbalpha_{k}=-\bbbeta_k$ for $k\geq0$, summing up the equations in \eqref{sorry1} and \eqref{sorry2} yields
\begin{equation}\label{sorry3}
(\bbA_s+\bbA_d)\bbx_{k+1}=2\bbz_{k+1}.
\end{equation}
Considering the definition of the oriented incidence matrix $\bbE_u=\bbA_s+\bbA_d$, we obtain that $\bbE_u \bbx_k=2\bbz_k$ holds for $k>0$. According to the initial condition $\bbE_u \bbx_0=2\bbz_0$, we can conclude that the relation $\bbE_u \bbx_k=2\bbz_k$ holds for $k\geq 0$.

Subtract the update for $\bbbeta_k$ in \eqref{sorry2} from the update for $\bbalpha_k$ in \eqref{sorry1} and consider the relation $\bbbeta_k=-\bbalpha_k$ to obtain
\begin{equation}\label{sorry4}
\bbalpha_{k+1}=\bbalpha_{k}+\frac{c}{2}(\bbA_s-\bbA_d)\bbx_{k+1}.
\end{equation}
Substituting $\bbA_s-\bbA_d$ in \eqref{sorry4} by $\bbE_o$ implies that
\begin{equation}\label{update_alpha}
\bbalpha_{k+1}=\bbalpha_{k}+\frac{c}{2}\bbE_o\bbx_{k+1}.
\end{equation}
 Hence, if $\bbalpha_k$ lies in the column space of matrix $\bbE_o$, then $\bbalpha_{k+1}$ also lies  in the column space of $\bbE_o$. According to the third condition of Assumption \ref{initial_val_assum}, $\bbalpha_0$ satisfies this condition, therefore $\bbalpha_k$ lies in the column space of matrix $\bbE_o$ for all $k\geq0$.


\section{Proof of Proposition \ref{update_system_prop}\label{app_update_system}}

The update for the multiplier $\bblambda$ in \eqref{ADMM_lambda_update} implies that we can substitute $\bblambda_k$ by $\bblambda_{k+1}-c( \bbA\bbx_{k+1}+\bbB\bbz_{k+1} )$ to simplify  \eqref{DQM_update_1} as
\begin{equation}\label{javad}
\nabla f(\bbx_{k})+
\bbH_{k}(\bbx_{k+1}-\bbx_{k})+\bbA^T\bblambda_{k+1}+c \bbA^{T}
\bbB\left(  \bbz_{k}-\bbz_{k+1}  \right)=\bb0.
\end{equation}
Considering the first result of Lemma \ref{lag_var_lemma} that $\bbalpha_k=-\bbbeta_k$ for $k\geq0$ in association with the definition $\bbA=[\bbA_s;\bbA_d]$ implies that the product $\bbA^T\bblambda_{k+1}$ is equivalent to
\begin{equation}\label{chand}
\bbA^T\bblambda_{k+1}=\bbA_s^T\bbalpha_{k+1}+\bbA_d^T\bbbeta_{k+1}=(\bbA_s-\bbA_d)^T\bbalpha_{k+1}.
\end{equation}
According to the definition $\bbE_o:=\bbA_s-\bbA_d$, the right hand side of \eqref{chand} can be simplified as
\begin{equation}\label{chand2}
\bbA^T\bblambda_{k+1}=\bbE_o^T\bbalpha_{k+1}.
\end{equation}
Based on the structures of the matrices $\bbA$ and $\bbB$, and the
definition $\bbE_u:=\bbA_s+\bbA_d$, we can simplify $\bbA^T\bbB$ as 
\begin{equation}\label{chand3}
\bbA^T\bbB=-\bbA_s^T-\bbA_d^T=-\bbE_u^T.
\end{equation}
Substituting the results in \eqref{chand2} and \eqref{chand3} into \eqref{javad} leads to
\begin{equation}\label{javad2}
\nabla f(\bbx_{k})+ \bbH_{k}(\bbx_{k+1}-\bbx_{k})+\bbE_o^T\bbalpha_{k+1}+c \bbE_u^{T}\left(  \bbz_{k+1}-\bbz_{k}  \right)=\bb0.
\end{equation}
The second result in Lemma \ref{lag_var_lemma} states that
$\bbz_k=\bbE_u\bbx_{k}/2$. Multiplying both sides of this equality
by $\bbE_u^{T}$ from left we obtain that $\bbE_u^{T}
\bbz_{k}=\bbE_u^{T}\bbE_u \bbx_{k}/2$ for $k\geq0$. Observing the
definition of the unoriented Laplacian
$\bbL_u:=\bbE_u^{T}\bbE_u/2$, we obtain that the product
$\bbE_u^{T} \bbz_{k}$ is equal to $\bbL_u\bbx_k$ for $k\geq0$.
Therefore, in \eqref{javad2} we can substitute $ \bbE_u^{T}\left(
\bbz_{k+1}-\bbz_{k}  \right)$ by $\bbL_u(\bbx_{k+1}-\bbx_k)$ and
write
\begin{equation}\label{javad3}
\nabla f(\bbx_{k})+ \left(\bbH_{k}+c\bbL_u\right)(\bbx_{k+1}-\bbx_{k})+\bbE_o^T\bbalpha_{k+1}=\bb0.
\end{equation}
Observe that the new variables $\bbphi_k$ are defined as $\bbphi_k:=\bbE_o^T\bbalpha_{k}$.
Multiplying both sides of \eqref{update_alpha} by $\bbE_o^T$ from the left hand side and considering the definition of oriented Laplacian $\bbL_o=\bbE_o^T\bbE_o/2$ follows the update rule of $\bbphi_{k}$ in \eqref{x_update_formula}, i.e.,
\begin{equation}\label{update_phi}
\bbphi_{k+1}=\bbphi_{k}+c\bbL_o\bbx_{k+1}.
\end{equation}
According to the definition $\bbphi_k=\bbE_o^T\bbalpha_{k}$ and the update formula in \eqref{update_phi}, we can conclude that $\bbE_o^T\bbalpha_{k+1}=\bbphi_{k+1}=\bbphi_{k}+c\bbL_o\bbx_{k+1}$. Substituting $\bbE_o^T\bbalpha_{k+1}$ by $\bbphi_{k}+c\bbL_o\bbx_{k+1}$ in \eqref{javad3} yields
\begin{equation}\label{javad4}
\nabla f(\bbx_{k})+ \left(\bbH_{k}+c\bbL_u\right)(\bbx_{k+1}-\bbx_{k})+\bbphi_{k}+c\bbL_o\bbx_{k+1}=\bb0.
\end{equation}
Observing the definition $\bbD=(\bbL_u+\bbL_o)/2$ we rewrite \eqref{javad4} as
\begin{equation}\label{javad5}
 \left(\bbH_{k}+2c\bbD\right)\bbx_{k+1}=\left(\bbH_{k}+c\bbL_u\right)\bbx_k-\nabla f(\bbx_{k})-\bbphi_{k}.
\end{equation}
Multiplying both sides of \eqref{javad5} by $ \left(\bbH_{k}+2c\bbD\right)^{-1}$ from the left hand side yields the first update in \eqref{x_update_formula}.


\section{Proof of Lemma  \ref{Hessian_Lipschitz_countinous}}\label{app_Hessian_Lipschitz}

Consider two arbitrary vectors $\bbx:=[\bbx_1;\dots;\bbx_n]  \in \reals^{np}$ and $\hbx:=[\hbx_1;\dots;\hbx_n] \in \reals^{np}$. Since the aggregate function Hessian is block diagonal where the $i$-th  diagonal block is given by $\nabla^2f_{i}(\bbx_i)$, we obtain that the difference of Hessians $\bbH(\bbx)-\bbH(\hbx)$ is also block diagonal where the $i$-th diagonal block ${\bbH(\bbx)_{ii}-\bbH(\hbx)}_{ii}$ is
\begin{equation}\label{matrix_difference}
{\bbH(\bbx)_{ii}-\bbH(\hbx)}_{ii}=
\nabla^2 f_{i}(\bbx_{i})-\nabla^2 f_{i}(\hbx_{i}) .
\end{equation}
Consider any vector $\bbv\in \reals^{np}$ and separate each $p$
components of vector $\bbv$ and consider it as a new vector called
$\bbv_{i}\in \reals^p$, i.e., $\bbv:=[\bbv_1;\dots;\bbv_n]$.
Observing the relation for the difference $\bbH(\bbx)-\bbH(\hbx)$
in \eqref{matrix_difference}, the symmetry of matrices
$\bbH(\bbx)$ and $\bbH(\hbx)$, and the definition of Euclidean
norm of a matrix that $\|\bbA\|=\sqrt{\lambda_{max}(\bbA^T\bbA)}$,
we obtain that the squared difference norm
$\|\bbH(\bbx)-\bbH(\hbx)\|^2$ can be written as
\begin{align}\label{inner_product_differnece}
\left\|\bbH(\bbx)-\bbH(\hbx)\right\|^2 \!
    &= \max_{\bbv} \frac{\bbv^T[\bbH(\bbx)-\bbH(\hbx)]^2\bbv}{\|\bbv\|^2} \\
    &= \max_{\bbv} \frac{\sum_{i=1}^n \bbv_{i}^T \left[\nabla^2 f_{i}(\bbx_{i})-\nabla^2 f_{i}(\hbx_{i})\right]^2    \bbv_{i}}                           {\|\bbv\|^2}\nonumber
\end{align}
Using the Cauchy-Schwarz inequality we can write
\begin{equation}\label{couchy_result}
\bbv_{i}^T\! \left[\nabla^2 f_{i}(\bbx_{i})\!-\!\nabla^2 f_{i}(\hbx_{i})\right]^2\! \bbv_{i} \leq
        \left\|\nabla^2 f_{i}(\bbx_{i})\!-\!\nabla^2 f_{i}(\hbx_{i})\right\|^2 \! \|\bbv_{i}\|^2
\end{equation}
Substituting the upper bound in \eqref{couchy_result} into \eqref{inner_product_differnece} implies that the squared norm $\left\|\bbH(\bbx)-\bbH(\hbx)\right\|^2$ is bounded above as
\begin{equation}\label{norm_difference_2}
\left\|\bbH(\bbx)-\bbH(\hbx)\right\|^2\leq
 \max_{\bbv} \frac{\sum_{i=1}^n\left\|\nabla^2 f_{i}(\bbx_{i})-\nabla^2 f_{i}(\hbx_{i})\right\|^2  \|\bbv_{i}\|^2}{\|\bbv\|^2}.
\end{equation}
Observe that Assumption 3 states that local objective functions Hessian $\nabla^2 f_{i}(\bbx_{i})$ are Lipschitz continuous with constant $L$, i.e. $\|\nabla^2 f_{i}(\bbx_{i})-\nabla^2 f_{i}(\hbx_{i})\|\leq L \|\bbx_{i}-\hbx_{i}\|$. Considering this inequality the upper bound in \eqref{norm_difference_2} can be changed by replacing $\|\nabla^2 f_{i}(\bbx_{i})-\nabla^2 f_{i}(\hbx_{i})\|$ by $L \|\bbx_{i}-\hbx_{i}\|$ which yields
\begin{equation}\label{alef}
\left\|\bbH(\bbx)-\bbH(\hbx)\right\|^2\leq
 \max_{\bbv} \frac{L^2\sum_{i=1}^n\left\|\bbx_{i}-\hbx_{i}\right\|^2  \|\bbv_{i}\|^2}{\sum_{i=1}^n\|\bbv_{i}\|^2}.
\end{equation}
Note that for any sequences of scalars such as $a_{i}$ and $b_i$, the inequality $\sum_{i=1}^n a_{i}^2b_i^2\leq (\sum_{i=1}^n a_{i}^2)(\sum_{i=1}^n b_{i}^2)$ holds. If we divide both sides of this relation by $\sum_{i=1}^n b_{i}^2$ and set $a_i=\|\bbx_{i}-\hbx_{i}\|$ and $b_i=\|\bbv_{i}\|$, we obtain
\begin{equation}\label{alef_2}
\frac{\sum_{i=1}^n\left\|\bbx_{i}-\hbx_{i}\right\|^2
\|\bbv_{i}\|^2}{\sum_{i=1}^n\|\bbv_{i}\|^2}
    \leq \sum_{i=1}^n\left\|\bbx_{i}-\hbx_{i}\right\|^2.
\end{equation}
Combining the two inequalities in \eqref{alef} and \eqref{alef_2} leads to
\begin{equation}\label{alef_3}
\left\|\bbH(\bbx)-\bbH(\hbx)\right\|^2
   \leq \max_{\bbv} L^2 \sum_{i=1}^n\left\|\bbx_{i}-\hbx_{i}\right\|^2  .
\end{equation}
Since the right hand side of \eqref{alef_3} does not depend on
$\bbv$ we can eliminate the maximization with respect to $\bbv$.
Further, note that according to the structure of vectors $\bbx$
and $\hbx$, we can write
$\left\|\bbx-\hbx\right\|^2=\sum_{i=1}^n\left\|\bbx_{i}-\hbx_{i}\right\|^2$.
These two observations in association with \eqref{alef_3} imply
that
\begin{equation}\label{rangarang}
\left\|\bbH(\bbx)-\bbH(\hbx)\right\|^2 \leq  L^2
\left\|\bbx-\hbx\right\|^2,
\end{equation}
Computing the square roots of terms in \eqref{rangarang} yields \eqref{H_Lipschitz_claim}.



\section{Proofs of Propositions \ref{error_vector_proposition} and \ref{error_vector_proposition_2}}\label{app_error_bound}

The fundamental theorem of calculus implies that the difference of gradients $\nabla f(\bbx_{k+1})-\nabla f(\bbx_k)$ can be written as
\begin{equation}\label{grad_dif}
\nabla f(\bbx_{k+1})-\nabla f(\bbx_k)=\int_{0}^1 \bbH(s\bbx_{k+1}+(1-s)\bbx_{k})(\bbx_{k+1}-\bbx_{k})\ ds.
\end{equation}
By computing norms of both sides of \eqref{grad_dif} and
considering that norm of integral is smaller than integral of norm
we obtain that 
\begin{equation}\label{grad_dif_norm}
\|\nabla f(\bbx_{k+1})-\nabla f(\bbx_k)\|\!\leq\!\!\int_{0}^1\!\!\!\!  \|\bbH(s\bbx_{k+1}+(1-s)\bbx_{k})(\bbx_{k+1}-\bbx_{k})\| ds.
\end{equation}
The upper bound $M$ for the eigenvalues of the Hessians as in \eqref{aggregate_hessian_eigenvlaue_bounds}, implies that $\| \bbH\left(s\bbx+(1-s)\hbx\right)(\bbx-\hbx)\|\leq M\|\bbx-\hbx\|$. Substituting this upper bound into \eqref{grad_dif_norm} leads to
\begin{equation}\label{grad_dif_norm_2}
\|\nabla f(\bbx_{k+1})-\nabla f(\bbx_k)\|\leq M \|\bbx_{k+1}-\bbx_{k}\|.
\end{equation}
The error vector norm $\|\bbe_{k}^{DLM}\|$ in \eqref{DLM_error} is bounded above as
 \begin{equation}\label{grad_dif_norm_dlm}
 \|\bbe_{k}^{DLM}\|\leq \|\nabla f(\bbx_{k+1})-\nabla f(\bbx_k)\|+\rho\|\bbx_{k+1}-\bbx_{k}\|.
 \end{equation}
By substituting the upper bound for $\|\nabla f(\bbx_{k+1})-\nabla f(\bbx_k)\|$ in \eqref{grad_dif_norm_2} into \eqref{grad_dif_norm_dlm}, the claim in \eqref{bound_for_DLM_error} follows.

 To prove \eqref{bound_for_DQM_error}, first we show that $ \|\bbe_{k}^{DQM}\|\leq2M\|\bbx_{k+1}-\bbx_{k}\|$ holds. Observe that the norm of error vector $\bbe_{k}^{DQM}$ defined \eqref{DQM_error} can be upper bounded using the triangle inequality as
 \begin{equation}\label{grad_dif_norm_dqm}
 \|\bbe_{k}^{DQM}\|\leq \|\nabla f(\bbx_{k+1})-\nabla f(\bbx_k)\|+\|\bbH_k(\bbx_{k+1}-\bbx_{k})\|.
 \end{equation}
Based on the Cauchy-Schwarz inequality and the upper bound $M$ for
the eigenvalues of Hessians as in
\eqref{aggregate_hessian_eigenvlaue_bounds}, we obtain $\|
\bbH_k(\bbx_{k+1}-\bbx_k)\|\leq M\|\bbx_{k+1}-\bbx_k\|$. Further,
as mentioned in \eqref{grad_dif_norm_2} the difference of
gradients $\|\nabla f(\bbx_{k+1})-\nabla f(\bbx_k)\|$ is upper
bounded by $M \|\bbx_{k+1}-\bbx_{k}\|$. Substituting these upper
bounds for the terms in the right hand side of
\eqref{grad_dif_norm_dqm} yields
 \begin{equation}\label{bound_1_dif_grad_dqm}
  \|\bbe_{k}^{DQM}\|\leq2M\|\bbx_{k+1}-\bbx_{k}\|.
 \end{equation}
The next step is to show that $ \|\bbe_{k}^{DQM}\|\leq(L/2)\|\bbx_{k+1}-\bbx_{k}\|^2$. Adding and subtracting the integral $\int_{0}^1\bbH(\bbx_k)(\bbx_{k+1}-\bbx_{k})\ ds$ to the right hand side of \eqref{grad_dif} results in
\begin{align}\label{grad_dif_new}
&\nabla f(\bbx_{k+1})-\nabla f(\bbx_k) =\int_{0}^1\bbH(\bbx_k)(\bbx_{k+1}-\bbx_{k})\ ds \nonumber \\
&+ \int_{0}^1 \left[\bbH(s\bbx_{k+1}+(1-s)\bbx_{k})-\bbH(\bbx_k)\right](\bbx_{k+1}-\bbx_{k})\ ds.
\end{align}
First observe that the integral $\int_{0}^1\bbH(\bbx_k)(\bbx_{k+1}-\bbx_{k})\ ds$ can be simplified as $\bbH(\bbx_k)(\bbx_{k+1}-\bbx_{k})$. Observing this simplification and regrouping the terms yield
\begin{align}\label{grad_dif_new2}
&\nabla f(\bbx_{k+1})-\nabla f(\bbx_k) -\bbH(\bbx_k)(\bbx_{k+1}-\bbx_{k})= \nonumber \\
& \int_{0}^1 \left[\bbH(s\bbx_{k+1}+(1-s)\bbx_{k})-\bbH(\bbx_k)\right](\bbx_{k+1}-\bbx_{k})\ ds.
\end{align}
 Computing norms of both sides of \eqref{grad_dif_new2}, considering the fact that norm of integral is smaller than integral of norm, and using Cauchy-Schwarz inequality lead to
\begin{align}\label{grad_dif_new3}
&\|\nabla f(\bbx_{k+1})-\nabla f(\bbx_k) -\bbH(\bbx_k)(\bbx_{k+1}-\bbx_{k})\|\leq \\
& \int_{0}^1 \left\|\bbH(s\bbx_{k+1}+(1-s)\bbx_{k})-\bbH(\bbx_k)\right\|\|\bbx_{k+1}-\bbx_{k}\| ds. \nonumber
\end{align}
Lipschitz continuity of the Hessian as in \eqref{H_Lipschitz_claim} implies that $\left\|\bbH(s\bbx_{k+1}+(1-s)\bbx_{k})-\bbH(\bbx_k)\right\|  \leq sL\|\bbx_{k+1}-\bbx_{k}\|$. By substituting this upper bound into the integral in \eqref{grad_dif_new3} and substituting the left hand side of \eqref{grad_dif_new3} by $\|\bbe_{k}^{DQM}\|$ we obtain
\begin{equation}\label{grad_dif_new4}
\|\bbe_{k}^{DQM}\|\leq \int_{0}^1 sL\|\bbx_{k+1}-\bbx_{k}\|^2 ds.
\end{equation}
Simplification of the integral in \eqref{grad_dif_new4} follows
\begin{equation}\label{grad_dif_new44}
\|\bbe_{k}^{DQM}\|\leq \frac{L}{2}\|\bbx_{k+1}-\bbx_{k}\|^2.
\end{equation}
The results in \eqref{bound_1_dif_grad_dqm} and \eqref{grad_dif_new44} follow the claim in \eqref{bound_for_DQM_error}.


\section{Proof of Lemma \ref{equalities_for_optima_lemma}\label{app_equalities_for_optima}}

In this section we first introduce an equivalent version of Lemma \ref{equalities_for_optima_lemma} for the DLM algorithm. Then, we show the validity of both lemmata in a general proof.

\begin{lemma}\label{equalities_for_optima_lemma_2}
Consider DLM as defined by \eqref{ADMM_z_update}-\eqref{DLM_x_update}. If Assumption \ref{initial_val_assum} holds true, then the optimal arguments $\bbx^*$, $\bbz^*$, and $\bbalpha^*$ satisfy
\begin{align}
   \nabla f(\bbx_{k+1})-\nabla f(\bbx^*)+\bbe_{k}^{DLM}
          + \bbE_o^T(\bbalpha_{k+1}-\bbalpha^*) \nonumber \\
         -c\bbE_u^T \left( \bbz_{k}-\bbz_{k+1}  \right)
   & = \bb0, \label{relation111} \\
   2(\bbalpha_{k+1}-\bbalpha_k)-{c}\bbE_o(\bbx_{k+1}-\bbx^*)
   &=\bb0, \label{relation222} \\
   \bbE_u(\bbx_{k}-\bbx^*)-2(\bbz_k-\bbz^*)
   &=\bb0.\label{relation333}
\end{align}
\end{lemma}
Notice that the claims in Lemmata \ref{equalities_for_optima_lemma} and  \ref{equalities_for_optima_lemma_2} are identical except in the error term of the first equalities. To provide a general framework to prove the claim in these lemmata we introduce $\bbe_k$ as the general error vector. By replacing $\bbe_k$ with $\bbe_k^{DQM}$ we obtain the result of DQM in Lemma \ref{equalities_for_optima_lemma} and by setting $\bbe_k=\bbe_k^{DLM}$ the result in Lemma  \ref{equalities_for_optima_lemma_2} follows. We start with the following Lemma that captures the KKT conditions of optimization problem \eqref{original_optimization_problem4}.


\begin{lemma}\label{KKT}
 Consider the optimization problem \eqref{original_optimization_problem4}. The optimal Lagrange multiplier $\bbalpha^*$, primal variable $\bbx^*$ and auxiliary variable $\bbz^*$ satisfy the following system of equations
\begin{align}\label{KKT_cliam}
\nabla f(\bbx^*)+\bbE_{o}^T\bbalpha^*=\bb0,\quad
\bbE_{o}\bbx^*=\bb0,\quad
\bbE_u\bbx^*=2\bbz^*.
\end{align}
\end{lemma}

\begin{myproof}
First observe that the KKT conditions of the decentralized optimization problem in \eqref{original_optimization_problem4} are given by
\begin{equation}\label{KKT1}
\nabla f(\bbx^*)+\bbA^T{\bblambda^*}=\bb0,\quad
\bbB^T{\bblambda^*}=\bb0,\quad
\bbA\bbx^*+\bbB\bbz^*=\bb0.
\end{equation}
Based on the definitions of the matrix
$\bbB=[-\bbI_{mp};-\bbI_{mp}]$ and the optimal Lagrange multiplier
$\bblambda^*:=[\bbalpha^*;\bbbeta^*]$, we obtain that
$\bbB^T{\bblambda^*}=\bb0$ in \eqref{KKT1} is equivalent to
$\bbalpha^*=-\bbbeta^*$. Considering this result and the
definition $\bbA=[\bbA_s;\bbA_d]$, we obtain
\begin{equation}\label{shepesh}
\bbA^T{\bblambda^*}=\bbA_s^T\bbalpha^*+\bbA_d^T\bbbeta^*=(\bbA_s-\bbA_d)^T\bbalpha^*.
\end{equation}
The definition $\bbE_o:=\bbA_s-\bbA_d$ implies that the right hand side of \eqref{shepesh} can be simplified as $\bbE_o^T\bbalpha^*$ which shows $\bbA^T{\bblambda^*}=\bbE_o^T\bbalpha^*$. Substituting $\bbA^T{\bblambda^*}$ by $\bbE_o^T\bbalpha^*$ into the first equality in \eqref{KKT1} follows the first claim in \eqref{KKT_cliam}.

Decompose the KKT condition $\bbA\bbx^*+\bbB\bbz^*=\bb0$ in \eqref{KKT1} based on the definitions of $\bbA$ and $\bbB$ as
\begin{equation}\label{fast1}
\bbA_s\bbx^*-\bbz=\bb0,\quad
\bbA_d\bbx^*-\bbz=\bb0.
\end{equation}
Subtracting the equalities in \eqref{fast1} implies that $(\bbA_s-\bbA_d)\bbx^*=\bb0$ which by considering the definition $\bbE_o=\bbA_s-\bbA_d$, the second equation in \eqref{KKT_cliam} follows. Summing up the equalities in \eqref{fast1} yields $(\bbA_s+\bbA_d)\bbx^*=2\bbz$. This observation in association with the definition $\bbE_u=\bbA_s-\bbA_d$ follows the third equation in \eqref{KKT_cliam}.
\end{myproof}


\textbf{Proofs of Lemmata \ref{equalities_for_optima_lemma} and \ref{equalities_for_optima_lemma_2}}:
First note that the results in Lemma \ref{lag_var_lemma} are also valid for DLM \cite{ling2014dlm}. Now, consider the first order optimality condition for primal updates of DQM and DLM in \eqref{DQM_update_1} and \eqref{DLM_optimality_cond}, respectively. Further, recall the definitions of error vectors $\bbe_k^{DQM}$ and $\bbe_k^{DLM}$ in \eqref{DQM_error} and \eqref{DLM_error}, respectively. Combining these observations we obtain that
\begin{equation}\label{whatever}
\nabla f(\bbx_{k+1})+\bbe_{k}+\bbA^T\bblambda_{k}+c \bbA^{T} \left( \bbA\bbx_{k+1}+\bbB \bbz_{k}  \right)=\bb0.
\end{equation}
Notice that by setting $\bbe_{k}=\bbe_{k}^{DQM}$ we obtain the
update for primal variable of DQM; likewise, setting
$\bbe_{k}=\bbe_{k}^{DLM}$ yields to the update of DLM.

Observe that the relation $\bblambda_{k}=\bblambda_{k+1}-c( \bbA\bbx_{k+1}+\bbB\bbz_{k+1} ) $ holds for both DLM and DQM according to to the update formula for Lagrange multiplier in \eqref{ADMM_lambda_update} and \eqref{DQM_update_2}. Substituting $\bblambda_k$ by $\bblambda_{k+1}-c( \bbA\bbx_{k+1}+\bbB\bbz_{k+1} ) $ in \eqref{whatever} follows
\begin{equation}\label{whatever2}
\nabla f(\bbx_{k+1})+\bbe_{k}+\bbA^T\bblambda_{k+1}+c \bbA^{T} \bbB\left(  \bbz_{k}-\bbz_{k+1}  \right)=\bb0
\end{equation}
Based on the result in Lemma \ref{lag_var_lemma}, the components of the Lagrange multiplier $\bblambda=[\bbalpha;\bbbeta]$ satisfy $\bbalpha_{k+1}=-\bbbeta_{k+1}$. Hence, the product $\bbA^T{\bblambda_{k+1}}$ can be simplified as $\bbA_s^T\bbalpha_{k+1}-\bbA_d^T\bbalpha_{k+1}=\bbE_o^T\bbalpha_{k+1}$ considering the definition that $\bbE_o=\bbA_s-\bbA_d$. Furthermore, note that according to the definitions we have that $\bbA=[\bbA_s;\bbA_d]$ and $\bbB=[-\bbI;-\bbI]$ which implies that $\bbA^T\bbB=-(\bbA_s+\bbA_d)^T=-\bbE_u^T$. By making these substitutions into \eqref{whatever2} we can write
\begin{align}\label{phone}
\nabla f(\bbx_{k+1})+\bbe_{k}+\bbE_o^T\bbalpha_{k+1}-c\bbE_u^T \left( \bbz_{k}-\bbz_{k+1}  \right)&=\bb0.
\end{align}
The first result in Lemma \ref{KKT} is equivalent to $\nabla f(\bbx^*)+\bbE_{o}^T\bbalpha^*=\bb0$. Subtracting both sides of this equation from the relation in \eqref{phone} follows the first claim of Lemmata \ref{equalities_for_optima_lemma} and \ref{equalities_for_optima_lemma_2}.

We proceed to prove the second and third claims in Lemmata \ref{equalities_for_optima_lemma} and \ref{equalities_for_optima_lemma_2}. The update formula for $\bbalpha_k$ in \eqref{update_alpha} and the second result in Lemma \ref{KKT}  that $\bbE_o\bbx^*=0$ imply that the second claim of Lemmata \ref{equalities_for_optima_lemma} and \ref{equalities_for_optima_lemma_2} are valid. Further, the result in Lemma \ref{lag_var_lemma} guaranteaes that $\bbE_u\bbx_{k}=2\bbz_k$. This result in conjunction with the result in Lemma \ref{KKT} that $\bbE_u\bbx^*=2\bbz^*$ leads to the third claim of Lemmata \ref{equalities_for_optima_lemma} and \ref{equalities_for_optima_lemma_2}.


\section{Proofs of Theorems \ref{DQM_convergence} and \ref{DLM_convergence}}\label{linear_convg_app}

To prove Theorems \ref{DQM_convergence} and \ref{DLM_convergence} we show a sufficient condition for the claims in these theorems. Then, we prove these theorems by showing validity of the sufficient condition. To do so, we use the general coefficient $\beta_k$ which is equivalent to $\zeta_k$ in the DQM algorithm and equivalent to $\rho+M$ in the DLM method. These definitions and the results in Propositions \ref{error_vector_proposition} and \ref{error_vector_proposition_2} imply that
\begin{equation}\label{new_realtion}
\|\bbe_k\|\leq \beta_k \|\bbx_{k+1}-\bbx_k\|,
\end{equation}
where $\bbe_k$ is $\bbe_k^{DQM}$ in DQM and $\bbe_k^{DLM}$ in DLM.
The sufficient condition of Theorems \ref{DQM_convergence} and \ref{DLM_convergence} is studied in the following lemma.

\begin{lemma}\label{equi_cond_linea_convg}
Consider the DLM and DQM algorithms as defined in \eqref{ADMM_z_update}-\eqref{DLM_x_update} and \eqref{DQM_update_1}-\eqref{DQM_update_2}, respectively. Further, conducer $\delta_k$ as a sequence of positive scalars. If Assumptions \ref{initial_val_assum}-\ref{Lipschitz_assumption} hold true then the sequence $\|\bbu_{k}-\bbu^*\|_\bbC^2$ converges linearly as
\begin{equation}\label{zarbe}
\|\bbu_{k+1}-\bbu^*\|_{\bbC}^2\  \leq\ \frac{1}{1+\delta_k}  \|\bbu_{k}-\bbu^*\|_{\bbC}^2,
\end{equation}
if the following inequality holds true,
\begin{align}\label{claim_2}
&\! \beta_k\|\bbx_{k+1}\!-\!\bbx^*\|\|\bbx_{k+1}\!-\!\bbx_k\|\!+\!\delta_k c \|\bbz_{k+1}\!-\!\bbz^*\|^2
\!+\!\frac{\delta_k}{c}\|\bbalpha_{k+1}\!-\!\bbalpha^*\|^2
\nonumber \\
&
\leq m\|\bbx_{k+1}\!-\!\bbx^*\|^2+c \|\bbz_{k+1}\!-\!\bbz_k\|^2+\frac{1}{c}\|\bbalpha_{k+1}\!-\!\bbalpha_k\|^2.
\end{align}

\end{lemma}

\begin{myproof}
Proving linear convergence of the sequence $\|\bbu_k-\bbu^*\|_\bbC^2$ as mentioned in \eqref{zarbe} is equivalent to showing that
 \begin{equation}\label{linear_convg}
\delta_k \|\bbu_{k+1}-\bbu^* \|_\bbC^2\leq  \|\bbu_{k}-\bbu^* \|_\bbC^2-\|\bbu_{k+1}-\bbu^* \|_\bbC^2.
 \end{equation}
 According to the definition $\|\bba\|_\bbC^2:=\bba^T\bbC\bba$ we can show that
 \begin{align}\label{mirror}
2(\bbu_{k}-\bbu_{k+1})^T\bbC(\bbu_{k+1}-\bbu^*)&= \|\bbu_{k}-\bbu^* \|_\bbC^2- \|\bbu_{k+1}-\bbu^*\|_\bbC^2 \nonumber\\
&\qquad - \|\bbu_{k}-\bbu_{k+1} \|_\bbC^2.
\end{align}
 The relation in \eqref{mirror} shows that the right hand side of \eqref{linear_convg} can be substituted by $2(\bbu_{k}-\bbu_{k+1})^T\bbC(\bbu_{k+1}-\bbu^*)+\|\bbu_{k}-\bbu_{k+1} \|_\bbC^2$. Applying this substitution into \eqref{linear_convg} leads to
 \begin{equation}\label{working}
 \delta_k \|\bbu_{k+1}-\bbu^* \|_\bbC^2 \leq 2(\bbu_{k}-\bbu_{k+1})^T\bbC(\bbu_{k+1}-\bbu^*)+ \|\bbu_{k}-\bbu_{k+1} \|_\bbC^2
 \end{equation}
This observation implies that to prove the linear convergence as claimed in \eqref{zarbe}, the inequality in \eqref{working} should be satisfied.

 We proceed by finding a lower bound for the term $2(\bbu_{k}-\bbu_{k+1})^T\bbC(\bbu_{k+1}-\bbu^*)$ in \eqref{working}. By regrouping the terms in \eqref{phone} and multiplying both sides of equality by $(\bbx_{k+1}-\bbx^*)^T$ from the left hand side we obtain that the inner product $(\bbx_{k+1}-\bbx^*)^T(\nabla f(\bbx_{k+1})-\nabla f(\bbx^*))$ is equivalent to
\begin{align}\label{glad}
&(\bbx_{k+1}-\bbx^*)^T(\nabla f(\bbx_{k+1})-\nabla f(\bbx^*))= \nonumber\\
&\qquad \qquad  -(\bbx_{k+1}-\bbx^*)^T\bbe_{k}-(\bbx_{k+1}-\bbx^*)^T\bbE_o^T(\bbalpha_{k+1}-\bbalpha^*)\nonumber\\
&\qquad \qquad+ c(\bbx_{k+1}-\bbx^*)^T\bbE_u^T ( \bbz_{k}-\bbz_{k+1}  ).
\end{align}
Based on \eqref{relation2}, we can substitute
$(\bbx_{k+1}-\bbx^*)^T\bbE_o^T(\bbalpha_{k+1}-\bbalpha^*)$ in
\eqref{glad} by
$(2/c)(\bbalpha_{k+1}-\bbalpha_{k})^T(\bbalpha_{k+1}-\bbalpha^*)$.
Further, the result in \eqref{relation3} implies that the term
$c(\bbx_{k+1}-\bbx^*)^T\bbE_u^T \left( \bbz_{k}-\bbz_{k+1}
\right)$ in \eqref{glad} is equivalent to $2c \left(
\bbz_{k}-\bbz_{k+1}\right)^T  (\bbz_{k+1}-\bbz^*)$. Applying these
substitutions into \eqref{glad} leads to 
\begin{align}\label{you}
&(\bbx_{k+1}\!-\!\bbx^*)^T(\nabla f(\bbx_{k+1})-\nabla f(\bbx^*))= -(\bbx_{k+1}\!-\!\bbx^*)^T\bbe_{k} \\
& \!\!+\frac{2}{c}(\bbalpha_{k}-\bbalpha_{k+1})^T(\bbalpha_{k+1}-\bbalpha^*)\!+\! 2c \left( \bbz_{k}-\bbz_{k+1}\right)^T\!(\bbz_{k+1}-\bbz^*).\nonumber
\end{align}
Based on the definitions of matrix $\bbC$ and vector $\bbu$ in \eqref{C_u_definitions}, the last two summands in the right hand side of \eqref{you} can be simplified as
\begin{align}\label{money}
&\frac{2}{c}(\bbalpha_{k}-\bbalpha_{k+1})^T(\bbalpha_{k+1}-\bbalpha^*)+2c \left( \bbz_{k}-\bbz_{k+1}\right)^T  (\bbz_{k+1}-\bbz^*)\nonumber\\
&\qquad =2(\bbu_{k}-\bbu_{k+1})^T\bbC(\bbu_{k+1}-\bbu^*).
\end{align}
Considering the simplification in \eqref{money} we can rewrite \eqref{you} as
\begin{align}\label{fine}
&(\bbx_{k+1}-\bbx^*)^T(\nabla f(\bbx_{k+1})-\nabla f(\bbx^*))
 \\
&\qquad \quad =-(\bbx_{k+1}-\bbx^*)^T\bbe_{k}+2(\bbu_{k}-\bbu_{k+1})^T\bbC(\bbu_{k+1}-\bbu^*). \nonumber
\end{align}
Observe that the objective function $f$ is strongly convex with
constant $m$ which implies the inequality
$m\|\bbx_{k+1}-\bbx^*\|^2\leq (\bbx_{k+1}-\bbx^*)^T(\nabla
f(\bbx_{k+1})-\nabla f(\bbx^*))$ holds true. Considering this
inequality from the strong convexity of objective function $f$ and
the simplification for the inner product
$(\bbx_{k+1}-\bbx^*)^T(\nabla f(\bbx_{k+1})-\nabla f(\bbx^*))$ in
\eqref{fine}, the following inequality holds 
\begin{equation}\label{strong_convexity_result}
m\|\bbx_{k+1}-\bbx^*\|^2 +(\bbx_{k+1}-\bbx^*)^T\bbe_{k}\leq2(\bbu_{k}-\bbu_{k+1})^T\bbC(\bbu_{k+1}-\bbu^*).
\end{equation}
Substituting the lower bound for the term
$2(\bbu_{k}-\bbu_{k+1})^T\bbC(\bbu_{k+1}-\bbu^*)$ in
\eqref{strong_convexity_result} into \eqref{working}, it follows that 
the following condition is sufficient to have \eqref{zarbe}, 
\begin{align}\label{new_claim_for_linear_convg}
 \delta_k \|\bbu_{k+1}-\bbu^* \|_\bbC^2& \leq m\|\bbx_{k+1}-\bbx^*\|^2 +(\bbx_{k+1}-\bbx^*)^T\bbe_{k}\nonumber \\ &\qquad + \|\bbu_{k}-\bbu_{k+1} \|_\bbC^2.
\end{align}
We emphasize that inequality \eqref{new_claim_for_linear_convg} implies the linear convergence result in \eqref{zarbe}. Therefore, our goal is to show that if  \eqref{claim_2} holds, the relation in \eqref{new_claim_for_linear_convg} is also valid and consequently the result in \eqref{zarbe} holds. According to the definitions of matrix $\bbC$ and vector $\bbu$ in \eqref{C_u_definitions}, we can substitute $ \|\bbu_{k+1}-\bbu^* \|_\bbC^2$ by $c \|\bbz_{k+1}-\bbz^*\|^2+({1}/{c})\|\bbalpha_{k+1}-\bbalpha^*\|^2$ and $\|\bbu_{k}-\bbu_{k+1} \|_\bbC^2$ by $c \|\bbz_{k+1}-\bbz_k\|^2+({1}/{c})\|\bbalpha_{k+1}-\bbalpha_k\|^2$. Making these substitutions into \eqref{new_claim_for_linear_convg} yields
\begin{align}\label{imp_inequality}
&\delta_k c \|\bbz_{k+1}-\bbz^*\|^2+\frac{\delta_k}{c}\|\bbalpha_{k+1}-\bbalpha^*\|^2  \leq m\|\bbx_{k+1}-\bbx^*\|^2   \\
&\qquad +(\bbx_{k+1}-\bbx^*)^T\bbe_{k}+c \|\bbz_{k+1}-\bbz_k\|^2+\frac{1}{c}\|\bbalpha_{k+1}-\bbalpha_k\|^2.\nonumber
\end{align}
%
The inequality in \eqref{new_realtion} implies that $-\|\bbe_k\|$ is lower bounded by $-\beta_{k}\|\bbx_{k+1}-\bbx_k\|$. This lower bound in conjunction with the fact that inner product of two vectors is not smaller than the negative of their norms product leads to
\begin{equation}\label{pen2}
(\bbx_{k+1}-\bbx^*)^T\!\bbe_{k}\geq-\beta_k\|\bbx_{k+1}-\bbx^*\|\|\bbx_{k+1}-\bbx_k\| .
\end{equation}
Substituting $(\bbx_{k+1}-\bbx^*)^T\bbe_{k}$ in \eqref{imp_inequality} by its lower bound in \eqref{pen2} leads to a sufficient condition for \eqref{imp_inequality} as in \eqref{claim_2}, i.e.,
\begin{align}\label{claim_2324}
&\! \beta_k\|\bbx_{k+1}\!-\!\bbx^*\|\|\bbx_{k+1}\!-\!\bbx_k\|\!+\!\delta_k c \|\bbz_{k+1}\!-\!\bbz^*\|^2
\!+\!\frac{\delta_k}{c}\|\bbalpha_{k+1}\!-\!\bbalpha^*\|^2
\nonumber \\
&
\leq m\|\bbx_{k+1}\!-\!\bbx^*\|^2+c \|\bbz_{k+1}\!-\!\bbz_k\|^2+\frac{1}{c}\|\bbalpha_{k+1}\!-\!\bbalpha_k\|^2.
\end{align}
Observe that if \eqref{claim_2324} holds true, then
\eqref{imp_inequality} and its equivalence
\eqref{new_claim_for_linear_convg} are valid and as a result the
inequality in \eqref{zarbe} is also satisfied. 
\end{myproof}

According to the result in Lemma \ref{equi_cond_linea_convg}, the
sequence $\|\bbu_{k}-\bbu^*\|^2$ converges linearly as mentioned
in \eqref{zarbe} if the inequality in \eqref{claim_2} holds true.
Therefore, in the following proof we show that for
\begin{equation}\label{general_delta}
\delta_k=\min
\Bigg\{
\frac{(\mu-1)(c\gamma_u^2-\eta_k\beta_k)\gamma_o^{2}}
     {\mu\mu'( c\Gamma_u^2\gamma_u^2+4\beta_k^2/c(\mu'-1))}
,
\frac{  m-{\beta_k}/{\eta_k} }{{ c}\Gamma_u^2/4
+\mu M^2/c\gamma_o^{2}}
\Bigg\},
\end{equation}
the inequality in \eqref{claim_2} holds and consequently \eqref{zarbe} is valid.

\textbf{Proofs of Theorems \ref{DQM_convergence} and \ref{DLM_convergence}}: we show that if the constant $\delta_k$ is chosen as in \eqref{general_delta}, then the inequality in \eqref{claim_2} holds true.
To do this first we should find an upper bound for
$\beta_k\|\bbx_{k+1}-\bbx^*\|\|\bbx_{k+1}-\bbx_k\| $ regarding the
terms in the right hand side of \eqref{claim_2}. Observing the
result of Lemma \ref{lag_var_lemma} that $\bbE_u\bbx_k=2\bbz_{k}$
for times $k$ and $k+1$, we can write 
\begin{equation}\label{salt}
\bbE_u(\bbx_{k+1}-\bbx_{k})=2(\bbz_{k+1}-\bbz_{k}).
\end{equation}
The singular values of $\bbE_u$ are bounded below by $\gamma_u$. Hence, equation \eqref{salt} implies that $\|\bbx_{k+1}-\bbx_{k}\|$ is upper bounded by
\begin{equation}\label{flower}
\|\bbx_{k+1}-\bbx_{k}\|\leq \frac{2}{\gamma_u} \|\bbz_{k+1}-\bbz_{k}\|.
\end{equation}
Multiplying both sides of \eqref{flower} by $\beta_k\|\bbx_{k+1}-\bbx^*\|$ yields
\begin{equation}\label{flower2}
\beta_k\|\bbx_{k+1}-\bbx^*\|\|\bbx_{k+1}-\bbx_{k}\|\leq \frac{2\beta_k}{\gamma_u}\|\bbx_{k+1}-\bbx^*\| \|\bbz_{k+1}-\bbz_{k}\|.
\end{equation}
Notice that for any vectors $\bba$ and $\bbb$ and positive constant $\eta_{k}>0$ the inequality $2\|\bba\|\|\bbb\|\leq (1/\eta_{k})\|\bba\|^2+\eta_{k}\|\bbb\|^2 $ holds true. By setting $\bba=\bbx_{k+1}-\bbx^*$ and $\bbb=(1/\gamma_u^2)(\bbz_{k+1}-\bbz_{k})$ the inequality $2\|\bba\|\|\bbb\|\leq (1/\eta_{k})\|\bba\|^2+\eta_{k}\|\bbb\|^2 $ is equivalent to
\begin{equation}\label{milad}
 \frac{2}{\gamma_u}\|\bbx_{k+1}-\bbx^*\| \|\bbz_{k+1}-\bbz_{k}\|\leq
 \frac{1}{\eta_{k}}\|\bbx_{k+1}-\bbx^*\|^2+\frac{\eta_{k} }{\gamma_u^2} \|\bbz_{k+1}-\bbz_{k}\|^2.
\end{equation}
Substituting the upper bound for $ ({2}/{\gamma_u})\|\bbx_{k+1}-\bbx^*\| \|\bbz_{k+1}-\bbz_{k}\|$ in \eqref{milad} into \eqref{flower2} yields
\begin{equation}\label{flower3}
\beta_k\|\bbx_{k+1}-\bbx^*\|\|\bbx_{k+1}-\bbx_{k}\|\!\leq \!
 \frac{\beta_k}{\eta_{k}}\|\bbx_{k+1}-\bbx^*\|^2+\frac{\eta_{k}\beta_k }{\gamma_u^2}
 \|\bbz_{k+1}-\bbz_{k}\|^2.
\end{equation}
Notice that inequality \eqref{flower3} provides an upper bound for
$\beta_k\|\bbx_{k+1}-\bbx^*\|\|\bbx_{k+1}-\bbx_{k}\|$ in
\eqref{claim_2} regarding the terms in the right hand side of
inequality which are $\|\bbx_{k+1}-\bbx^*\|^2$ and
$\|\bbz_{k+1}-\bbz_{k}\|^2$. The next step is to find upper bounds
for the other two terms in the left hand side of \eqref{claim_2}
regarding the terms in the right hand side of \eqref{claim_2}
which are $\|\bbx_{k+1}-\bbx^*\|^2$, $\|\bbz_{k+1}-\bbz_{k}\|^2$,
and $\|\bbalpha_{k+1}-\bbalpha_k\|^2$. First we start with $
\|\bbz_{k+1}-\bbz^*\|^2$. The relation in \eqref{relation3} and
the upper bound $\Gamma_u$ for the singular values of matrix
$\bbE_u$ yield
\begin{equation}\label{second_term_bound}
\delta_k c \|\bbz_{k+1}-\bbz^*\|^2 \leq \frac{\delta_k c\Gamma_u^2}{4} \|\bbx_{k+1}-\bbx^*\|^2.
\end{equation}
The next step is to bound $(\delta_k/c)\|\bbalpha_{k+1}-\bbalpha^*\|$ in terms of the term in the right hand side of (47). First, note that for any vector $\bba$, $\bbb$, and $\bbc$, and constants $\mu$ and $\mu'$ which are larger than 1, i.e. $\mu,\mu'>1$, we can write
\begin{align}\label{kahn}
(1-\frac{1}{\mu'})(1-\frac{1}{\mu})\|\bbc\|^2
&\leq
\|\bba+\bbb+\bbc\|^2+(\mu'-1)\|\bba\|^2\nonumber\\
&\qquad + (\mu-1)(1-\frac{1}{\mu'})\|\bbb\|^2 .
\end{align}
Set $\bba=c\bbE_u^T ( \bbz_{k}-\bbz_{k+1})$, $\bbb=\nabla f(\bbx^*)-\nabla f(\bbx_{k+1})$, and $\bbc=\bbE_o^T(\bbalpha^*-\bbalpha_{k+1})$. By choosing these values and observing equality \eqref{relation1} we obtain $\bba+\bbb+\bbc=\bbe_{k}$. Hence, by making these substitutions for $\bba$, $\bbb$, $\bbc$, and $\bba+\bbb+\bbc$ into \eqref{kahn} we can write
\begin{align}\label{beats}
(1-\frac{1}{\mu'})(1-\frac{1}{\mu})&\|\bbE_o^T(\bbalpha_{k+1}-\bbalpha^*)\|^2
\leq  \|\bbe_{k}\|^2
\\
&+(\mu'-1)\|c\bbE_u^T ( \bbz_{k}-\bbz_{k+1})\|^2 \nonumber\\
&+(\mu-1)(1-\frac{1}{\mu'})\|\nabla f(\bbx_{k+1})-\nabla f(\bbx^*)\|^2 \nonumber.
\end{align}
Notice that according to the result in Lemma \ref{lag_var_lemma}, the Lagrange multiplier $\bbalpha_{k}$ lies in the column space of $\bbE_o$ for all $k\geq0$. Further, recall that the optimal multiplier $\bbalpha^*$ also lies in the column space of $\bbE_o$. These observations show that $\bbalpha^*-\bbalpha_k$ is in the column space of $\bbE_o$. Hence, there exits a vector $\bbr\in \reals^{np}$ such that $\bbalpha^*-\bbalpha_k=\bbE_o\bbr$. This relation implies that $\|\bbE_o^T(\bbalpha_{k+1}-\bbalpha^*)\|^2$ can be written as $\|\bbE_o^T\bbE_o\bbr\|^2=\bbr^T(\bbE_o^T\bbE_o)^2\bbr$. Observe that since the eigenvalues of matrix $(\bbE_o^T\bbE_o)^2$ are the squared of eigenvalues of the matrix $\bbE_o^T\bbE_o$, we can write $\bbr^T(\bbE_o^T\bbE_o)^2\bbr \geq \gamma_o^2\bbr^T\bbE_o^T\bbE_o\bbr$, where $\gamma_o$ is the smallest non-zero singular value of the oriented incidence matrix $\bbE_o$. 
Observing this inequality and the definition $\bbalpha^*-\bbalpha_k=\bbE_o\bbr$ we can write
\begin{equation}\label{paper100}
\left\|\bbE_o^T(\bbalpha_{k+1}-\bbalpha^*)\right\|^2 \geq \gamma_o^2 \| \bbalpha_{k+1}-\bbalpha^* \|^2.
\end{equation}
Observe that the error norm $\|\bbe_{k}\|$ is bounded above by $\beta_{k}\|\bbx_{k+1}-\bbx_k\|$ as in \eqref{new_realtion} and the norm $\|c\bbE_u^T ( \bbz_{k}-\bbz_{k+1})\|^2$ is upper bounded by $ c^2\Gamma_u^2\|\bbz_{k}-\bbz_{k+1}\|^2$ since all the singular values of the unoriented matrix $\bbE_u$ are smaller than $\Gamma_u$. Substituting these upper bounds and the lower bound in \eqref{paper100} into \eqref{beats} implies 
\begin{align}\label{coffee}
& (1-\frac{1}{\mu'})(1-\frac{1}{\mu})\gamma_o^2\|\bbalpha_{k+1}-\bbalpha^*\|^2
\leq   \beta_k^2  \|\bbx_{k+1}-\bbx_k\|^2  \\
&+\!(\mu'-1)c^2\Gamma_u^2\| \bbz_{k}-\bbz_{k+1}\|^2 \!+\!(\mu-1)(1-\frac{1}{\mu'})M^2\|\bbx_{k+1}\!-\!\bbx^*\|^2 \nonumber
\end{align}
Considering the result in \eqref{flower}, $\|\bbx_{k+1}-\bbx_{k}\|$ is upper by $ ({2}/{\gamma_u} )\|\bbz_{k+1}-\bbz_{k}\|$. Therefore, we can substitute $\|\bbx_{k+1}-\bbx_{k}\|$ in the right hand side of \eqref{coffee} by its upper bound $ ({2}/{\gamma_u} )\|\bbz_{k+1}-\bbz_{k}\|$. Making this substitution, dividing both sides  by $ (1-{1}/{\mu'})(1-{1}/{\mu})\gamma_o^2$, and regrouping the terms lead to
\begin{align}\label{third_term_bound}
\|\bbalpha_{k+1}-\bbalpha^*\|^2& \leq
\frac{\mu M^2}{\gamma_o^2}\|\bbx_{k+1}-\bbx^*\|^2
\\
&\!\!\!\!\!\!\!\!\!\!\!
+ \left[\frac{4\mu\mu'\beta_k^2}{\gamma_u^2\gamma_o^2(\mu-1)(\mu'-1)}\!+\!\frac{\mu\mu'c^2\Gamma_u^2}{(\mu-1)\gamma_o^2}\right] \!\! \| \bbz_{k}-\bbz_{k+1}\|^2.\nonumber
\end{align}
Considering the upper bounds for $\beta_k\|\bbx_{k+1}-\bbx^*\|\|\bbx_{k+1}-\bbx_k\|$, $\|\bbz_{k+1}-\bbz^*\|^2$, and $\|\bbalpha_{k+1}-\bbalpha_{k}\|^2$, in \eqref{flower3}, \eqref{second_term_bound}, and \eqref{third_term_bound}, respectively, we obtain that if the inequality
\begin{align}\label{final}
&
\left[
\frac{\beta_k}{\eta_{k}}
+\frac{\delta_k c\Gamma_u^2}{4}
+\frac{\delta_k\mu M^2}{c\gamma_o^2}
\right]
\|\bbx_{k+1}-\bbx^*\|^2
+ \\
&
\left[
\frac{4\delta_k\mu\mu'\beta_k^2}{c\gamma_u^2\gamma_o^2(\mu-1)(\mu'-1)} +\frac{\delta_k\mu\mu'c\Gamma_u^2}{(\mu-1)\gamma_o^2}
+\frac{\eta_{k}\beta_k}{\gamma_u^2}
\right]
\|\bbz_{k+1}-\bbz_k\|^2\nonumber\\
&
\qquad \leq m\|\bbx_{k+1}-\bbx^*\|^2+c \|\bbz_{k+1}-\bbz_k\|^2+\frac{1}{c}\|\bbalpha_{k+1}-\bbalpha_k\|^2.\nonumber
\end{align}
holds true, \eqref{claim_2} is satisfied. Hence, the last step is to show that for the specific choice of $\delta_k$ in \eqref{general_delta} the result in \eqref{final} is satisfied. In order to make sure that \eqref{final} holds, it is sufficient to show that the coefficients of $\|\bbx_{k+1}-\bbx^*\|^2$ and $\|\bbz_{k+1}-\bbz_k\|^2$ in the left hand side of \eqref{final} are smaller than the ones in the right hand side. Hence, we should verify the validity of inequalities
\begin{align}
&\qquad\qquad  \qquad \frac{\beta_k}{\eta_{k}}
+\frac{\delta_k c\Gamma_u^2}{4}
+\frac{\delta_k\mu M^2}{c\gamma_o^2}
 \leq
  m,  \label{hat1}\\
&\frac{4\delta_k\mu\mu'\beta_k^2}{c\gamma_u^2\gamma_o^2(\mu-1)(\mu'-1)} +\frac{\delta_k\mu\mu'c\Gamma_u^2}{(\mu-1)\gamma_o^2}
+\frac{\eta_{k}\beta_k}{\gamma_u^2}
 \leq
c.\label{hat2}
\end{align}
Considering the inequality for $\delta_k$ in \eqref{general_delta} we obtain that \eqref{hat1} and \eqref{hat2} are satisfied. Hence, if $\delta_k$ satisfies condition in \eqref{general_delta}, \eqref{final} and consequently \eqref{claim_2} are satisfied. Now recalling the result of Lemma \ref{equi_cond_linea_convg} that inequality \eqref{claim_2} is a sufficient condition for the linear convergence in \eqref{zarbe}, we obtain that the linear convergence holds. By setting $\beta_k=\zeta_k$ we obtain the linear convergence of DQM in Theorem \ref{DQM_convergence} is valid and the linear coefficient in \eqref{general_delta} can be simplified as \eqref{DLM_DQM_delta}. Moreover, setting $\beta_k=\rho+M$ follows the linear convergence of DLM as in Theorem \ref{DLM_convergence} with the linear constant in \eqref{DLM_delta}.

\end{appendices}